\documentclass[11pt,a4paper,reqno]{amsart}
\usepackage{amsmath,amssymb,mathtools,mathrsfs}
\usepackage[T1]{fontenc}
\usepackage[utf8]{inputenc}
\usepackage{lmodern}
\usepackage[textwidth=162mm,textheight=237mm,centering,headsep=6mm]{geometry}
\usepackage{microtype}
\usepackage{longtable,booktabs,array,calc}
\usepackage{fvextra,needspace,etoolbox}
\usepackage{xurl}
\usepackage[unicode,hidelinks]{hyperref}
\usepackage{bookmark}
\AtBeginEnvironment{longtable}{\footnotesize}
\DefineVerbatimEnvironment{verbatim}{Verbatim}{breaklines=true,breakanywhere=true,fontsize=\footnotesize}

\makeatletter
\patchcmd\longtable{\par}{\if@noskipsec\mbox{}\fi\par}{}{}
\newcommand{\ResultHeading}[1]{\if@noskipsec\leavevmode\par\fi\par\addvspace{5pt}\noindent\textbf{#1}}
\newcommand{\ProofHeading}[1]{\needspace{4\baselineskip}\if@noskipsec\leavevmode\par\fi\par\addvspace{3pt}\noindent\textit{#1}}
\newcommand{\ResultLabel}[2]{\def\@currentlabel{#1}\label{#2}}
\makeatother
\title[Fourteenth-canonical birationality]{Fourteenth-canonical birationality of threefolds with positive geometric genus}
\author{Zhengjie Yu}
\date{10 November 2026}
\hypersetup{pdftitle={Fourteenth-canonical birationality of threefolds with positive geometric genus},pdfauthor={Zhengjie Yu},pdflang={en}}
\begin{document}
\begin{abstract}
We prove that the pluricanonical map of a complex projective threefold of general type with canonical singularities and positive geometric genus is birational in every degree at least fourteen. The bound is sharp. 
\end{abstract}
\maketitle
\markboth{Zhengjie Yu}{Fourteenth-canonical birationality}
{\small\tableofcontents}
\medskip\normalfont
\needspace{4\baselineskip}

\hypertarget{introduction}{%
\needspace{10\baselineskip}\hypertarget{section-1}{}\section{1. Introduction}\label{introduction}}

Pluricanonical maps provide a natural way to study varieties of general type through their canonical rings. A basic problem is to determine a degree, depending only on the dimension, beyond which these maps are birational. For smooth surfaces of general type, Bombieri proved that every \(m\)-canonical map is birational for \(m\ge5\), and this bound is optimal \cite{ref3}. In arbitrary dimension, the work of Hacon--McKernan, Takayama, and Tsuji establishes the existence of a uniform birationality bound \cite{HM2006,Takayama2006,Tsuji1999}. Determining sharp explicit bounds remains a separate problem.

In dimension three, singularities of minimal models and the possible scarcity of low-degree pluricanonical sections play an essential role. Reid's Riemann--Roch formula expresses their contribution in terms of baskets of terminal singularities \cite{ref25}. Chen--Chen developed this into a systematic method for controlling plurigenera and canonical volumes, together with restriction arguments on surfaces and curves \cite{ref9,ref8,ref7}. These methods connect the numerical possibilities for a canonical ring with the geometry of the first moving pluricanonical systems.

We study the sharp birationality bound under the assumption of positive geometric genus. The case of geometric genus one requires particular care: the canonical system consists of a single divisor, and the first moving system may occur in degree four.

We work over the complex numbers. For a projective threefold \(X\) with canonical singularities, put
\[
P_m(X)=h^0(X,\omega_X^{[m]}),\qquad p_g(X)=P_1(X),
\]
and let \(\varphi_m\) be the rational map defined by the complete \(m\)-canonical system. Here \(\omega_X^{[m]}=(\omega_X^{\otimes m})^{**}\). Our main result is the following.

\needspace{3\baselineskip}

\hypertarget{result-1.1}{}\ResultHeading{Theorem 1.1.} Let \(X\) be a complex projective threefold of general type with canonical singularities. If \(p_g(X)>0\), then \(\varphi_m\) is birational onto its image for every integer \(m\ge14\). The integer fourteen is optimal.

A general hypersurface \(X_{28}\subset\mathbb P(1,3,4,5,14)\) has \(p_g=1\) and a nonbirational thirteenth-canonical map; see \hyperlink{section-7}{Section 7}. Since a nonzero canonical section embeds the ratios defining \(\varphi_m\) into those defining \(\varphi_{m+1}\), it suffices to prove birationality in degree fourteen.

The bounds improve substantially when low-degree plurigenera are larger. Chen proved birationality in degree five for \(p_g\ge4\), in degree six for \(p_g=3\), and in degree eight for \(p_g=2\) \cite[Theorem 1.2 and Section 4]{Chen2003}. In particular, degree eight is a sharp uniform bound for \(p_g\ge2\); see \cite[Example 1.2(2)]{ref14}. Chen--Chen proved that \(P_2\ge2\) implies birationality in every degree at least eleven \cite[Theorem 1.8]{ref7}.

For \(p_g=1\), Chen--Chen obtained birationality in every degree at least eighteen \cite[Corollary 1.7]{ref7}. Chen--Hu--Penegini improved the bound to seventeen and showed that failure in degree sixteen is confined to five possible baskets \cite[Theorem 1.3]{ref14}. Motivated by the weighted hypersurface above, they conjectured that degree fourteen suffices \cite[Conjecture B]{ref14}. Theorem~1.1 proves this conjecture.

Several important cases were already known. Hu--Yan proved fourteenth-canonical birationality under the assumption \(P_3\ge2\), and showed that failure in degree thirteen forces a pencil of \((1,2)\) surfaces \cite[Theorem 1.1]{ref19}. The irregular case admits a stronger result: the fifth-canonical map is birational \cite[Theorem 1.3]{ref11}. Chen--Zuo proved the degree-fourteen assertion when \(\chi(\mathcal O_X)\le0\) \cite[Theorem 1.1]{ref12}. Since a nonzero canonical section gives \(P_1\le P_2\le P_3\), the remaining case is
\[
P_1=P_2=P_3=1,\qquad q=0,\qquad\chi(\mathcal O_X)>0.
\]
On a smooth model \(W\), these equalities give \(h^{2,0}(W)=\chi(\mathcal O_W)>0\), so the residual case supplies nonzero holomorphic two-forms.

The use of two-forms in the study of pluricanonical maps has an earlier precedent. Luo studied their zero loci under birational transformations and, under conditions on those loci, obtained restrictions on terminal singularities and pluricanonical systems \cite{Luo2000}. Here we construct a global polar class from a two-form and a pluricanonical section. Its nonvanishing gives a plurigenus inequality and detects a reduced component of the canonical divisor; these provide the additional information needed after the numerical and restriction-theoretic reductions.

Our proof combines canonical restriction with holomorphic two-forms.
In Section~3.3, we associate to \(0\ne a\in H^0(mK)\) a natural quadratic map
\[
H^0(W,\Omega_W^2)\longrightarrow H^0((2m+2)K)/a^2H^0(2K)
\]
whose only zero is the origin. The construction and its nonvanishing are proved there.
It yields \(P_{2m+2}-P_2\ge h^{2,0}(W)\). When \(P_1=P_3=1\),
its residue detects a reduced canonical component. This gives
\(P_{2r}>P_{2r-1}\) and excludes powers of the first new fourth-canonical
section becoming divisible by the canonical generator.

Choose a canonical form \(s\) and a holomorphic two-form \(\alpha\), and
define \(v\) by \(\iota_vs=\alpha\). The normalized spaces
\(E_m=s^{-m}H^0(mK)\) satisfy \(v(E_m)\subset E_{m+2}\).
Mixed Hodge theory annihilates the action of \(v\) on every
finite-dimensional \(v\)-stable space. The saturated quotient by \(s\),
with its full multiplicities retained, relates this derivation to
the polar residue. Canonical restriction reduces a possible failure
to a four-canonical pencil of \((1,2)\) surfaces and the sixth map to
a quadratic map with rational image. A marked Weierstrass point on
the canonical genus-two curve controls the poles of normalized sections.
The global nonvanishing and the seventh restricted system exclude
one plurigenus configuration. In the remaining configuration,
an alternating Jacobian lies outside the quadratic invariant field;
its pole calculation produces a fourteenth-canonical section.

In the local calculations, we expand meromorphic functions already defined on the threefold. Their convergent Laurent expansions allow us to compare the first possible nonintegral exponents. Thus the argument requires neither a classification of differential equations nor an existence theorem for formal solutions. \hyperlink{section-2}{Section 2} recalls restriction and marked-curve estimates. \hyperlink{section-3}{Section 3} develops the differential and graded-ring tools. \hyperlink{section-4}{Section 4} reduces the problem to the marked \((1,2)\) pencil; Sections 5 and 6 prove \hyperlink{result-1.1}{Theorem 1.1}. \hyperlink{section-7}{Section 7} discusses sharpness. The appendix records the finite basket calculation needed at one numerical endpoint.

\smallskip
\noindent\textbf{AI assistance.}
During this work, the author used ChatGPT (OpenAI) to assist with checking arguments, symbolic and exact arithmetic computations, and language editing. The author began thinking about this research topic one year ago. The author takes full responsibility for the mathematical content, references, and final presentation.

\newcommand{\AcknowledgementsText}{The author is grateful to his adviser Meng Chen for his great encouragement and support. He is also grateful for his fellow students for their helpful discussion and friendship.
}
\ifdefempty{\AcknowledgementsText}{}{%
  \smallskip\noindent\textbf{Acknowledgements.} \AcknowledgementsText\par}

\needspace{4\baselineskip}

\hypertarget{preliminaries}{%
\needspace{10\baselineskip}\hypertarget{section-2}{}\section{2. Preliminaries}\label{preliminaries}}

\needspace{4\baselineskip}

\hypertarget{models-canonical-rings-and-baskets}{%
\hypertarget{section-2.1}{}\subsection{2.1. Models, canonical rings, and baskets}\label{models-canonical-rings-and-baskets}}

We work over the complex numbers. A fibration has connected fibres. We use the standard definitions of terminal and canonical singularities in \cite{ref22}. The discrepancy \(a(E,X)\) is the coefficient of \(E\) in \(K_{X'}-\pi^*K_X\); its log discrepancy is \(1+a(E,X)\).

Let \(X\) be a projective threefold of general type with canonical singularities. Its pluricanonical ring is unchanged on a resolution. The minimal model programme and finite generation \cite{ref2,ref22} provide a terminal \(\mathbb Q\)-factorial minimal model \(Y\) and a canonical model \(Z\). We use \(Y\) for nef intersections, \(Z\) for the canonical ring, and a smooth projective model \(W\) for differential forms. On a common resolution \(\pi:W\to Y\), set
\[
Q=\pi^*K_Y,\qquad V=Q^3=\operatorname{vol}(K_X),\qquad
h=h^0(W,\Omega_W^2),\qquad q=h^1(W,\mathcal O_W).
\]
Thus \(Q\) is nef and big and \(\chi(\mathcal O_W)=1-q+h-p_g\). Canonical divisors are chosen compatibly, so that \(K_W=Q+E_\pi\) with \(E_\pi\ge0\).

Write \(\omega_Z^{[m]}=(\omega_Z^{\otimes m})^{**}\), with the usual dual interpretation for negative \(m\), and put
\[
R=\bigoplus_{m\ge0}R_m,\qquad
R_m=H^0(Z,\omega_Z^{[m]})=H^0(W,mK_W),\qquad P_m=\dim R_m.
\]
For a fixed nonzero \(s\in R_1\), let
\[
E_m=R_m/s^m\subset\mathbb C(W),\qquad L_m=\mathbb C(E_m).
\]
Here \(E_m\) is a vector space of rational functions, whereas \(L_m\) is the field it generates. Multiplication by \(s\) gives \(E_m\subset E_{m+1}\). Moreover, \(\varphi_m\) is birational if and only if \(L_m=\mathbb C(W)\). In particular, when \(p_g>0\), birationality in degree \(m\) implies birationality in every larger degree.

We will need to distinguish inequalities between actual divisors from inequalities up to rational linear equivalence. The notation \(D\ge_{\mathbb Q}D'\) means that \(D-D'\) is \(\mathbb Q\)-linearly equivalent to an effective rational divisor. An inequality without this subscript compares coefficients. Round-up always applies to an actual representative. Neither rounding nor a coefficientwise comparison is inferred from numerical equivalence alone.

A terminal singularity has a basket of cyclic quotient singularities. We use the convention
\[
(b,r)=\frac1r(1,-1,b),\qquad 0<b\le r/2,\qquad (b,r)=1.
\]
Multiplicities are retained. Reid's Riemann--Roch formula and vanishing give, for \(m\ge2\),
\[
P_m=\frac{m(m-1)(2m-1)}{12}V-(2m-1)\chi
 +\sum_{(b,r)\in B}\sum_{j=1}^{m-1}
 \frac{\overline{jb}(r-\overline{jb})}{2r},                 \tag{2.1}\label{eq:2.1}
\]
where \(\overline{jb}\) is the least nonnegative residue modulo \(r\) \cite{ref25,ref9}. Consequently
\[
V=2P_2+6\chi-\sum_{(b,r)\in B}\frac{b(r-b)}r.              \tag{2.2}\label{eq:2.2}
\]
A packing replaces \((b,r),(c,u)\) by the primitive pair \((b+c,r+u)\). At fixed \(P_2,\chi\), it decreases the volume by
\[
\frac{b^2}{r}+\frac{c^2}{u}-\frac{(b+c)^2}{r+u}
=\frac{(bu-cr)^2}{ru(r+u)}.
\]
It is prime when \(|bu-cr|=1\). We use the canonical unpacking procedure of \cite[Section 3]{ref9}. A basket constituent need not itself be an actual singular point; assertions about a local Cartier index refer to one actual point on \(Y\).

\needspace{4\baselineskip}

\hypertarget{restriction-to-surfaces-and-curves}{%
\hypertarget{section-2.2}{}\subsection{2.2. Restriction to surfaces and curves}\label{restriction-to-surfaces-and-curves}}

For a smooth moving surface \(F\subset W\) of general type, let \(\sigma:F\to S\) be its smooth minimal-model contraction and put
\[
D=Q|_F,\qquad A=\sigma^*K_S,\qquad d=D^2.
\]
A surface of type \((a,b,c)\) has \(K_S^2=a\), \(p_g(S)=b\), and \(q(S)=c\).

\needspace{3\baselineskip}

\hypertarget{result-2.1}{}\ResultHeading{Proposition 2.1 (canonical restriction).} Suppose that \(F\) is semiample and \(Q\ge_{\mathbb Q}\mu F\), where \(\mu>0\). Then
\[
D\ge_{\mathbb Q}\frac{\mu}{1+\mu}A,\qquad V\ge\mu d.
                                                               \tag{2.3}\label{eq:2.3}
\]
If \(Q-\mu F\) is only pseudoeffective, the first inequality holds with every positive rational \(p<\mu\) in place of \(\mu\), and
\[
V\ge\frac{\mu^3}{(1+\mu)^2}K_S^2.
\]

\ProofHeading{Proof.} The first assertion is the canonical restriction inequality of \cite[Corollary 2.3]{ref10}, applied with the semiample Weil divisor \(K_Y\) and coefficient \(1/\mu\). We obtain the second by intersecting \(Q-\mu F\) with \(Q^2\). To treat the pseudoeffective case, we write
\[
Q-pF=\frac p\mu(Q-\mu F)+\left(1-\frac p\mu\right)Q
\]
is big, hence \(\mathbb Q\)-effective. Nefness of \(D,A\) gives \(d\ge p^2A^2/(1+p)^2\). Let \(p\) tend to \(\mu\). \(\square\)

We next recall the model of a pencil constructed by Chen--Hu--Jiang. Its usefulness here is that adjunction gives a comparison between actual effective divisors. We keep track of the model, since this comparison is used later before further blowups.

\needspace{3\baselineskip}

\hypertarget{result-2.2}{}\ResultHeading{Proposition 2.2 (a model of a pencil).} A rational pencil on \(Y\), with
general member \(F_Y\), admits a terminal \(\mathbb Q\)-factorial model
\(\mu:W_{\mathrm{sp}}\to Y\) resolving it, with smooth general fibre \(F\),
such that
\[
G=\mu^*(K_Y+F_Y)-K_{W_{\mathrm{sp}}}-F\ge0.\tag{2.4}\label{eq:2.4}
\]
If \(aK_Y\sim bF_Y+Z_Y\) is an integral relation with \(Z_Y\ge0\), then
\[
(1+a/b)D\sim_{\mathbb Q}K_F+G|_F+b^{-1}\mu^*Z_Y|_F.\tag{2.5}\label{eq:2.5}
\]
Equality \((1+a/b)D^2=DA\) implies \(V=(b/a)D^2\).

For the complete connected four-canonical pencil, choose the actual
\(K_Y=\operatorname{div}_Y(s)\). The corresponding actual divisors satisfy
\[
0\le D\le K_F\le5D,\qquad5D\ge\sigma^*A_0,\quad
A_0=\sigma_*K_F\in|K_S|.\tag{2.6}\label{eq:2.6}
\]

\ProofHeading{Proof.} The existence of the model and \eqref{eq:2.4} are proved in \cite[Proposition 3.1]{ref13}. By adjunction we obtain \eqref{eq:2.5}. If equality holds, its effective remainder and \(K_F-A\) are
\(D\)-orthogonal, so \(K_YF_YZ_Y=0\). Put \(z=K_Y^2Z_Y\) and \(v=K_Y^3\).
The divisor relation gives \(K_YZ_Y^2=az\), while \(0\le z<av\).
Hodge index gives \(z^2\ge avz\), hence \(z=0\) and \(aV=bD^2\).

For \eqref{eq:2.6}, write \(4\mu^*K_Y=F_s+E_4\), where \(F_s\) is the full fibre
defined by \(s^4\) and \(E_4\) is the fixed part. If \(Z_Y\) is the fixed
divisor downstairs, then
\(\mu^*F_Y=F+E_4-\mu^*Z_Y\) and \(E_4|_F=4D\). Equation \eqref{eq:2.4} therefore
gives \(5D-K_F=G|_F+\mu^*Z_Y|_F\ge0\); terminality gives \(K_F-D\ge0\).
Finally \(K_F=\sigma^*A_0+E_\sigma\). \(\square\)

Further blowups preserve the pulled-back inequalities and the
intersections \(D^2,DA\). In \eqref{eq:2.6} one pulls back \(K_F\) as well; one does
not replace it by the canonical divisor of the new fibre.

\needspace{3\baselineskip}

\hypertarget{result-2.3}{}\ResultHeading{Proposition 2.3 (curve restriction).} Suppose \(Q\ge_{\mathbb Q}\mu F\), and let \(C\) be a general irreducible member of a free moving curve system on \(F\). If \(D\ge_{\mathbb Q}\beta C\), set
\[
\xi=DC,\qquad \alpha(m)=\left(m-1-\mu^{-1}-\beta^{-1}\right)\xi.
\]
For \(\alpha(m)>1\),
\[
m\xi\ge2g(C)-2+\lceil\alpha(m)\rceil.                     \tag{2.7}\label{eq:2.7}
\]
If the \(m\)-canonical system separates general surfaces \(F\) and its restriction separates general curves \(C\), then \(\alpha(m)>2\) implies birationality. In particular,
\[
\xi\ge\frac{2g(C)-2}{1+\mu^{-1}+\beta^{-1}}.
\]
In the four-canonical pencil of \hyperlink{result-2.2}{Proposition 2.2}, for \(m\ge6\),
\[
|mK_W||_F\succeq|K_F+\lceil(m-5)D\rceil|\succeq|K_F|.    \tag{2.8}\label{eq:2.8}
\]
Thus the outer pencil and any canonical pencil on \(F\) are separated in degree \(m\ge6\).

\ProofHeading{Proof.} The first assertions are \cite[Theorem 2.11 and Remark 2.12]{ref8}; the formulation for sub-linear systems is also recorded in \cite[Theorem 2.1]{ref14}. The numerical lower bound follows by taking \(m\) large in \eqref{eq:2.7}. For \eqref{eq:2.8}, resolve the fractional support and write \(4Q=F_s+E_4\). Put \(L=(m-1)Q-E_4=(m-5)Q+F_s\). Then \(L-F\equiv(m-5)Q\) is nef and big. Kawamata--Viehweg vanishing gives restriction for \(K_W+\lceil L\rceil\), whose restriction is the first system on the right in \eqref{eq:2.8}. The inclusion upstairs follows from \(E_4\ge0\) and \(K_W+\lceil(m-1)Q\rceil\le mK_W\). Effectivity of \(D\) gives the second inclusion. Powers of \(s\) retain the four-pencil in every degree at least four. \(\square\)

\needspace{3\baselineskip}

\hypertarget{result-2.4}{}\ResultHeading{Lemma 2.4 (rounding).} Let \(x>0\) be rational and let \(c>n\ge1\)
be integers. If
\(mx\ge n+\lceil(m-c)x\rceil\) for all sufficiently large integers \(m\)
and \(\gcd(c,n+1)>1\), then \(x\ge(n+1)/c\).

\ProofHeading{Proof.} We write \(x=a/b\) in lowest terms. For \(b>1\), choose arbitrarily
large \(m\) with \((m-c)a\equiv1\pmod b\). The inequality gives
\(ca\ge(n+1)b-1\). The difference is divisible by \(\gcd(c,n+1)>1\),
so it is nonnegative. If \(b=1\), then \(x\ge1\ge(n+1)/c\). \(\square\)

We also use the following form of Langer's separation theorem \cite[Theorem 4.1]{ref23}. If \(N\) is a nef and big rational divisor on a smooth surface, \(N^2>8\), and \(|K_F+\lceil N\rceil|\) is not birational, there is a curve \(C\) through a very general point such that
\[
NC\le\frac4{1+\sqrt{1-8/N^2}}.                           \tag{2.9}\label{eq:2.9}
\]
Indeed, the cluster invariant for two distinct smooth points is eight. Choose the parameter \(b<\sqrt{N^2/8}\) in that theorem and let it approach the endpoint. Intersections with a rational divisor are discrete, so strict inequality opposite to \eqref{eq:2.9} would separate two very general points. This statement concerns the actual divisor \(N\); it does not compare \(N\) and its round-up on a negative curve.

On a genus-two curve, the pole semigroup at a Weierstrass point \(p\) is \(\langle2,5\rangle\). Thus functions with only a pole at \(p\), of order at most four, lie in the hyperelliptic field. A birational system with pole divisor bounded by \(mp\) requires \(m\ge5\). By Riemann--Roch the gaps are one and three; the least nonconstant pole has order two, and the first function outside the degree-two map field has order five.

\needspace{3\baselineskip}

\hypertarget{result-2.5}{}\ResultHeading{Lemma 2.5 (a marked canonical curve).} Suppose a fibration on \(W\)
has general fibre whose minimal model is a \((1,2)\) surface. Resolve its
canonical pencil and let \(C\) be a general member. For the actual divisor
\(Q=\pi^*\operatorname{div}_Y(s)\) there is a Weierstrass point \(p\in C\)
such that, with \(\xi=QC\),
\[
Q|_C=\xi p,\qquad
(a/s^m)|_C\in H^0(C,\mathcal O_C(\lfloor m\xi\rfloor p))
\quad(a\in R_m).\tag{2.10}\label{eq:2.10}
\]
Consequently, birationality of \(\varphi_m\) implies \(m\xi\ge5\).

\ProofHeading{Proof.} We write \(K_S=M+Z\) on the minimal surface. Hodge index and
adjunction parity give \(K_SM=M^2=1\); the canonical-null fixed part \(Z\)
has square zero and is therefore zero. The pencil has one simple base
point. After blowing it up, the actual canonical divisor restricts to
\(2p\) on a general free member \(C\), with \(K_C\sim2p\). Further point
blowups are avoided by a general member. Since \(Q\) and \(K_W-Q\) restrict
effectively, their sum \(2p\) gives the first equality. Pulling back the
effective divisor \(\operatorname{div}(a/s^m)+mK_Y\) gives the pole bound.
The Weierstrass semigroup \(\langle2,5\rangle\) proves the last assertion.
\(\square\)

\needspace{3\baselineskip}

\hypertarget{result-2.6}{}\ResultHeading{Lemma 2.6 (the marked local index).} Let \(P=\pi(p)\) be the image of the marked point of \hyperlink{result-2.5}{Lemma 2.5}, and let \(r_P\) be its local terminal index. Then \(r_P\xi\in\mathbb{Z}\), and \(r_P\) occurs as a basket denominator. If \(P\) is a cyclic quotient \(\frac1{r_P}(1,-1,b)\) and \(\xi=a/e\) is reduced, the canonical coefficient on its smooth index-one cover has multiplicity at most \(a\) at the origin. In particular, for \(r_P\ge3\) and \(\xi=1/3\) it has a linear term, which requires \(b=1\) in the normalized convention.

\ProofHeading{Proof.} The Cartier divisor \(r_PQ\) has integral restriction coefficient at \(p\), proving integrality. The terminal-type table \cite[Remark 2.1]{ref6} records the actual local index among its basket denominators: a \(cA/r\) point contributes its axial weight many copies of \((b,r)\), a \(cAx/4\) point contributes an index-four entry, and the remaining indices are two or three.

On a smooth index-one cover write \(s=f\omega_\#\), where \(\omega_\#=dx\wedge dy\wedge dz\) has character \(b\), so \(f\) has the faithful character \(-b\). If \(\tau\) is the local frame of \(\omega^{[r_P]}\) induced by \(\omega_\#^{r_P}\), put \(h=s^{r_P}/\tau\). The cover is the normalization of \(T^{r_P}=h\). By the actual divisor equality \eqref{eq:2.10}, its pullback to the marked disc is \(h=z^{r_Pa/e}\cdot\text{unit}\). Its normalized branches satisfy \(z=t^e\), \(f=t^a\cdot\text{unit}\) and lift to arcs centred at the cover origin. This forces \(\operatorname{mult}_0f\le a\). For \(a=1\), character \(-b\) must occur among the three linear characters \(1,-1,b\). Normalization and coprimality give \(b=1\). This criterion is unchanged by the inverse convention for the quotient parameter. \(\square\)

\needspace{4\baselineskip}

\hypertarget{holomorphic-two-forms-and-canonical-rings}{%
\needspace{10\baselineskip}\hypertarget{section-3}{}\section{3. Holomorphic two-forms and canonical rings}\label{holomorphic-two-forms-and-canonical-rings}}

\needspace{4\baselineskip}

\hypertarget{holomorphic-two-forms}{%
\hypertarget{section-3.1}{}\subsection{3.1. Holomorphic two-forms}\label{holomorphic-two-forms}}

A mixed Hodge structure has a weight filtration \(W_\bullet\). Cohomology of a smooth projective variety is pure of its cohomological degree. Morphisms, including cup products, are strict for the weight filtration. If \(U\) is the complement of a divisor in a smooth projective variety, \(H^1(U)\) has weights one and two, and its weight-one part is the image of \(H^1(W)\) \cite[Theorem 2.3.5 and Corollaries 3.2.13(ii), 3.2.14, and 3.2.17]{ref17}. We apply these assertions to the full rational cohomology groups and then complexify.

\needspace{3\baselineskip}

\hypertarget{result-3.1}{}\ResultHeading{Lemma 3.1 (closed rational forms).} Let \(W\) be smooth projective with \(q(W)=0\).

\begin{enumerate}
\def\labelenumi{(\alph{enumi})}
\item
  If \(\eta,\zeta\) are closed rational one-forms and \(\eta\wedge\zeta\) extends holomorphically, that extension is zero.
\item
  If \(\alpha_i\) are holomorphic two-forms, \(\eta_i\) are closed rational one-forms, and \(\sum_i\eta_i\wedge\alpha_i\) extends to a holomorphic three-form, that extension is zero.
\end{enumerate}

\ProofHeading{Proof.} We remove the poles and resolve the boundary, obtaining \(U\) with a smooth projective compactification and simple normal crossing boundary. Since \(H^1(W)=0\), \(H^1(U)\) is pure of weight two. The cup-product maps
\[
H^1(U)\otimes H^1(U)\longrightarrow H^2(U),\qquad
H^1(U)\otimes H^2(W)\longrightarrow H^3(U)
\]
have pure weight-four sources. By strictness their images meet \(W_2H^2(U)\), respectively \(W_3H^3(U)\), only in zero. The indicated holomorphic forms come from projective cohomology and belong to these lower weight pieces. Their classes therefore vanish.

To pass from cohomology classes to forms, we use logarithmic Hodge-to-de Rham degeneration \cite[Corollaries 3.2.13--3.2.14]{ref17}. After resolving the boundary, a global holomorphic \(p\)-form is logarithmic, and the map from global logarithmic \(p\)-forms to \(H^p(U,\mathbb C)\) is injective. Thus a global holomorphic form whose class vanishes on \(U\) is zero. This proves both assertions. The rational forms need not have rational periods or logarithmic poles. \(\square\)

We now return to threefolds. Fix nonzero \(s\in H^0(W,K_W)\) and \(\alpha\in H^0(W,\Omega_W^2)\), and define \(v\) by \(\iota_vs=\alpha\). The identities we shall use are
\[
df\wedge\alpha=v(f)s,\qquad L_vs=d\alpha=0.              \tag{3.1}\label{eq:3.1}
\]

\needspace{3\baselineskip}

\hypertarget{result-3.2}{}\ResultHeading{Lemma 3.2 (functions and their derivatives).} Assume \(q(W)=0\).

\begin{enumerate}
\def\labelenumi{(\alph{enumi})}
\item
  If \(v(f)\ne0\), then \(f,v(f)\) are algebraically independent.
\item
  Every finite-dimensional complex vector subspace of \(\mathbb C(W)\) stable under \(v\) is annihilated by \(v\).
\item
  If \(f:W\to\mathbb P^1\) has regular general fibre \(F\), restriction of holomorphic two-forms is injective and \(h\le p_g(F)-1\). If instead \(p_g(F)=1\), all global two-forms restrict to zero, and their generic evaluation rank is at most \(q(F)\).
\end{enumerate}

\ProofHeading{Proof.} If \(df\wedge d(v(f))=0\) and \(v(f)\ne0\), then \(\eta=df/v(f)\) is closed and \(\eta\wedge\alpha=s\), contradicting \hyperlink{result-3.1}{Lemma 3.1}. The characteristic-zero differential criterion gives (a).

For (b), a nonzero eigenvalue of \(v\) on a finite-dimensional stable space would contradict (a). If a nonzero nilpotent Jordan block remains, we can choose \(v(g)=h\ne0\), \(v(h)=0\), and hence \(v(g/h)=1\), again impossible. This proves (b). Notice that the conclusion allows nonconstant first integrals.

For (c), restrictions of global closed forms are flat cohomology sections on the smooth locus of the family. A two-form vanishing on one general fibre vanishes on all smooth fibres. The relative differential sequence writes it locally as \(dt\wedge\eta\), where \(\eta|_F\) is a holomorphic one-form. Thus \(q(F)=0\) gives injectivity. Restrictions of a basis remain independent over \(\mathbb C(t)\), since a relation at a general fibre would give a constant relation among flat sections. Write \(s=dt\wedge s_F\). If \(s_F\) belonged to their \(\mathbb C(t)\)-span, then \(s=\sum_i a_i(t)dt\wedge\alpha_i\), contradicting \hyperlink{result-3.1}{Lemma 3.1}. Hence \(h\le p_g(F)-1\).

If \(p_g(F)=1\), each \(dt\wedge\alpha\) is a base-rational multiple of \(s\). A nonzero multiple contradicts \hyperlink{result-3.1}{Lemma 3.1}, so all restrictions vanish. Dividing by \(dt\) now gives relative holomorphic one-forms, of generic rank at most \(q(F)\). \(\square\)

\needspace{3\baselineskip}

\hypertarget{the-canonical-quotient-and-differentiation}{%
\hypertarget{section-3.2}{}\subsection{3.2. The canonical quotient and differentiation}\label{the-canonical-quotient-and-differentiation}}

Let \(\mathcal A=R/(s)\), and denote by \(\mathcal S\) the ring of regular functions on \(\operatorname{Spec}\mathcal A\) minus its vertex. This is its graded saturation. It retains nilpotents: the canonical divisor need not be reduced.

\needspace{3\baselineskip}

\hypertarget{result-3.4}{}\ResultHeading{Proposition 3.3 (the canonical quotient).} Assume \(q=0\), and put \(h=h^{2,0}(W)\). Then \(\mathcal A\to\mathcal S\) is injective and
\[
(\mathcal S/\mathcal A)_m=
\begin{cases}\mathbb C^h&m=2,\\0&m\ne2.\end{cases}            \tag{3.2}\label{eq:3.3}
\]
For \(\alpha\ne0\), there is a degree-two derivation \(\Delta:R\to R\), with \(\Delta(s)=0\), and an element \(\theta\in\mathcal S_2\) such that
\[
\Delta_m(b)=s^{m+2}v(b/s^m),\qquad
\overline{\Delta_m(b)}=-m\theta\bar b.                    \tag{3.3}\label{eq:3.4}
\]
If additionally \(P_1=P_2=P_3=1\) and \(P_4=2\), choose \(x\) independent of \(s^4\) and put \(t=x/s^4\). Then \(\theta^2=\lambda\bar x\), for some \(\lambda\in\mathbb C\), and
\[
v^2(b/s^m)-m(m+2)\lambda t(b/s^m)\in E_{m+3}.             \tag{3.4}\label{eq:3.5}
\]

\ProofHeading{Proof.} All reflexive powers \(\mathcal O_Z(mK_Z)\) are Cohen--Macaulay. Locally, take the index-one cover. It is Gorenstein canonical and hence Cohen--Macaulay, and its eigensheaves are the reflexive canonical powers. A system of parameters is regular on every direct summand; twisting by an invertible index multiple handles all integers \(m\). Weil-divisor Kawamata--Viehweg vanishing and Serre duality \cite[Theorems 2.70 and 5.71]{ref22} give
\[
H^1(Z,\mathcal O_Z(mK_Z))=
\begin{cases}\mathbb C^h&m=1,\\0&m\ne1.\end{cases}
\]
For \(m\ge2\) use vanishing, for \(m<0\) use its dual form, and for \(m=0,1\) use rational singularities, \(q=0\), and \(h=h^{2,0}(W)\).

The full section ring of an integral Weil divisor on a normal variety is normal. Indeed, after choosing a rational trivialization, it is the intersection of the polynomial ring over the function field with the weighted Gauss valuation rings imposing all divisorial section inequalities. Hence \(R\) is \(S_2\), and its quotient by the nonzero divisor \(s\) is \(S_1\). Thus \(H^0_{\mathfrak m}(\mathcal A)=0\). The punctured-spectrum sequence gives
\[
0\longrightarrow\mathcal A\longrightarrow\mathcal S
\longrightarrow H^1_{\mathfrak m}(\mathcal A)\longrightarrow0.
\]
The local/sheaf cohomology comparison identifies
\(H^2_{\mathfrak m}(R)_m=H^1(Z,\mathcal O_Z(mK_Z))\).
For a non-standard grading, apply the usual comparison to a sufficiently divisible Veronese and each residue-class module. Applying local cohomology to
\(0\to R(-1)\xrightarrow{s}R\to\mathcal A\to0\) proves \eqref{eq:3.3}.

In a local canonical frame \(\nu\), write \(s=g\nu\) and \(b=f\nu^m\). Direct differentiation gives
\[
\Delta_m(b)=\nu^{m+1}\otimes
\bigl((g\,df-mf\,dg)\wedge\alpha\bigr).                   \tag{3.5}\label{eq:3.6}
\]
The first-order expression \(g\,df-mf\,dg\) is the weighted Gaussian map \(\gamma_{m,1}\) of Ballico--Fontanari \cite[formula (1) and Lemma 1]{BF2006}, wedged with \(\alpha\). This is holomorphic and independent of the frame: the terms containing its logarithmic derivative cancel. It therefore belongs to \(R_{m+2}\); the product rule proves that \(\Delta\) is a derivation.

For \(a\in R_p,b\in R_q\), the same formula gives
\[
qb\Delta_p(a)-pa\Delta_q(b)\in sR_{p+q+1}.
\]
Consequently \(-\Delta(\bar a)/(p\bar a)\) agrees on the positive-degree principal opens of the punctured quotient cone. These fractions define \(\theta\in\mathcal S_2\) and prove \eqref{eq:3.4}. Locally its representative modulo \(g\) is \((dg\wedge\alpha)\otimes\nu\).

Under the extra plurigenus assumptions,
\(\mathcal S_4=\mathcal A_4=\mathbb C\bar x\), giving
\(\theta^2=\lambda\bar x\). Applying \eqref{eq:3.4} twice yields
\[
\Delta^2(b)-m(m+2)\lambda xb\in sR_{m+3}.
\]
Division by \(s^{m+4}\) proves \eqref{eq:3.5}. \(\square\)

\needspace{3\baselineskip}

\hypertarget{result-3.5}{}\ResultHeading{Lemma 3.4 (generic divisorial test).} At a smooth generic point of
\(T\subset\operatorname{div}_Z(s)\), let
\(\ell=\operatorname{ord}_T(s)>0\). We call \(T\) a transverse component if
\(\alpha|_T\ne0\). If \(b\in R_m\) has order \(c\ge0\), then
\[
\operatorname{ord}_T\Delta_m(b)\ge\ell+c-1.
\]
At a tangent component this improves to \(\ell+c\); at a transverse
component equality holds when \(c-m\ell\ne0\). In the generic Artin ring
of the full divisor \(s=0\), \(\theta\) is zero at a tangent
component and has exact order \(\ell-1\) at a transverse component.
A nonzero product \(\theta\bar b\) is therefore detected at a transverse
component where \(b\) is a unit. If \(\theta^2=0\), each transverse
component on \(Z\) has multiplicity at least two.

\ProofHeading{Proof.} We write \(g=\rho^\ell h\), with \(h\) a unit, and
\(\iota_w\nu=\alpha\). The coefficient in \eqref{eq:3.6} is
\(gw(f)-mf\,w(g)\). If \(w(\rho)\) is a unit, its leading term is a
nonzero unit times \(c-m\ell\). If \(w(\rho)\) is divisible by \(\rho\),
differentiation loses no normal order. These are the transverse and
tangent cases. The same calculation of \(w(g)\bmod g\) proves the
assertions about \(\theta\). Generic localizations detect zero because
\(\mathcal A\) is \(S_1\); away from the vertex, its saturation has the
same local rings. \(\square\)

\needspace{4\baselineskip}

\hypertarget{polar-squares-and-reduced-canonical-components}{%
\hypertarget{section-polar}{}\subsection{3.3. Polar squares and reduced canonical components}\label{polar-squares-and-reduced-canonical-components}}
We call a prime divisor \(T\) visible to a two-form \(\alpha\) if its pullback to a smooth model of \(T\) is nonzero. At a smooth generic point this is precisely the transverse case above. We construct a global polar class by gluing a quadratic differential expression modulo \(a^2\). Its local polynomial is formally a normalized second Rankin--Cohen self-bracket of weight \(2m\) \cite[Section~1]{Zagier1994}. The following proposition gives the construction and proves its nonvanishing.

\hypertarget{result-polar-inequality}{}\ResultHeading{Proposition 3.5 (polar inequality).}\ResultLabel{3.5}{prop:polar-inequality}
Let \(Y\) be a projective terminal threefold with \(K_Y\) nef and big, and let \(W\) be a smooth projective resolution. Put \(R_j=H^0(Y,\mathcal O_Y(jK_Y))\), \(P_j=\dim R_j\), and \(H=H^0(W,\Omega_W^2)\). For every \(m\ge1\) and \(0\ne a\in R_m\), there is a natural homogeneous quadratic map
\[
\mathfrak P_a:H\longrightarrow R_{2m+2}/a^2R_2
\]
whose only zero is the origin. In particular, \(P_{2m+2}-P_2\ge\dim H\).

\ProofHeading{Proof.}
Write \(\mathcal F_j=\mathcal O_Y(jK_Y)\) and \(\mathcal Q_a=\mathcal F_{2m+2}/a^2\mathcal F_2\). On the smooth locus choose a canonical frame \(\nu\), and write \(a=f\nu^m\), \(\alpha=\iota_\delta\nu\). Set
\[
F_m(f,\delta)=(\delta f)^2-\frac{2m}{2m+1}f\delta^2f.
\]
If \(\nu'=g\nu\), then \(f'=g^{-m}f\), \(\delta'=g^{-1}\delta\), and
\[
\begin{aligned}
F_m(f',\delta')={}&g^{-2m-2}F_m(f,\delta)\\
&+\frac{m^2f^2}{2m+1}\bigl(2g^{-2m-3}\delta^2g-3g^{-2m-4}(\delta g)^2\bigr).
\end{aligned}
\]
Thus \(F_m(f,\delta)\nu^{2m+2}\) defines a section of \(\mathcal Q_a\) on \(Y_{\rm reg}\). The canonical powers are maximal Cohen--Macaulay: on an index-one cover they are direct summands of the finite pushforward of its Cohen--Macaulay structure sheaf. The exact sequence
\[
0\longrightarrow\mathcal F_2\xrightarrow{\ a^2\ }\mathcal F_{2m+2}\longrightarrow\mathcal Q_a\longrightarrow0
\]
gives depth at least two for \(\mathcal Q_a\) at every singular point. Terminal threefold singularities are isolated, so the section extends uniquely over them. Weil-divisor Kawamata--Viehweg vanishing and Serre duality give \(H^1(Y,\mathcal F_2)=0\) \cite[Theorem~2.70]{ref22}. Hence its global section space is \(R_{2m+2}/a^2R_2\), and the construction gives \(\mathfrak P_a\), quadratically in \(\alpha\).

We prove nonvanishing. Choose the resolution to be an isomorphism over \(Y_{\rm reg}\). The pullback of \(K_Y\) is big. A Kodaira decomposition as a positive ample class plus an effective class shows that its pairing with the closed semipositive form \(\alpha\wedge\overline\alpha\) is strictly positive for \(\alpha\ne0\). Exceptional divisors are invisible: functorial pullback of reflexive forms factors their restriction through their zero-dimensional centres \cite{ref21}. Therefore the effective divisor \((a)\sim mK_Y\) contains a visible prime \(T\). At its generic point write \(f=t^ku\), where \(u\) and \(\delta t\) are units. The leading coefficient of \(F_m\), in order \(2k-2\), is
\[
\left(k^2-\frac{2m}{2m+1}k(k-1)\right)u^2(\delta t)^2
=\frac{k(k+2m)}{2m+1}u^2(\delta t)^2\ne0.
\]
Thus every lift \(u_a(\alpha)\in R_{2m+2}\) has order \(2k-2\) at \(T\), whereas every member of \(a^2R_2\) has order at least \(2k\). At an invisible component of multiplicity \(k\), \(\delta\) preserves its prime ideal, and the same expression has order at least \(2k\). In particular \(\mathfrak P_a(\alpha)\ne0\).

If \(H\ne0\), the resulting system of quadrics has no base point on \(\mathbb P(H)\). Its morphism to \(\mathbb P(R_{2m+2}/a^2R_2)\) pulls back \(\mathcal O(1)\) to \(\mathcal O(2)\), so its fibres are finite. The target has dimension at least \(\dim H-1\). If \(H=0\), multiplication by \(a^2\) already gives the inequality.
\(\square\)

A visible prime on \(Y\) is not contracted by the canonical morphism: its image would have dimension at most one, so its two-form restriction would vanish by the same functorial pullback. Thus the generic residue calculations agree on the minimal and canonical models.

\hypertarget{result-canonical-polar}{}\ResultHeading{Corollary 3.6 (canonical polar nonvanishing).}\ResultLabel{3.6}{cor:canonical-polar}
Assume \(q=0\), \(P_1=P_3=1\), and \(0\ne\alpha\in H^0(W,\Omega_W^2)\). Let \(s\) generate \(R_1\), and let \(\theta\in\mathcal S_2\) be its Euler residue.
\begin{enumerate}
\def\labelenumi{(\alph{enumi})}
\item The divisor of \(s\) has a reduced visible prime \(T\). The residue \(\theta\) is a unit at its generic point. Thus \(\theta^r\ne0\) for every \(r\ge1\), and \(P_{2r}>P_{2r-1}\) for \(r\ge2\).
\item If \(P_4=2\), then \(h=1\). For any \(R_4=\mathbb Cs^4\oplus\mathbb Cx\), one has \(\theta^2=\lambda\bar x\) with \(\lambda\ne0\), and \(x\) is a unit on \(T\). In particular,
\[
x^j\notin sR_{4j-1}\quad(j\ge1),\qquad v_\alpha(x/s^4)\notin E_5.
\]
Every visible component on which \(x\) is a unit is reduced.
\end{enumerate}

\ProofHeading{Proof.}
Let \(u_s\in R_4\) lift \(\mathfrak P_s(\alpha)\). If all visible components of \((s)\) had multiplicity \(k\ge2\), the valuation calculation in Proposition~\hyperlink{result-polar-inequality}{3.5} would give order at least \(k\) for \(u_s\) at every component of \((s)\); invisible components have order at least \(2k\). At every other prime, \(s\) is a unit, so \(u_s/s\) has no pole there. Normal reflexive extension would then give \(u_s/s\in R_3\). Since \(P_3=1\), this would imply \(u_s\in\mathbb Cs^4\subset s^2R_2\), a contradiction. Hence a reduced visible component \(T\) exists. Locally \(\theta\) is \((\delta f)\nu^2\bmod f\), where \(s=f\nu\), so it is a unit at \(T\). Its square is the class of \(u_s\) modulo \(s\). All its powers are nonzero. \hyperlink{result-3.4}{Proposition 3.3} identifies \(\mathcal S_{2r}\) with \(\mathcal A_{2r}=R_{2r}/sR_{2r-1}\) for \(r\ge2\), proving (a).

Multiplication by \(s\) gives \(P_2=1\). The polar inequality gives \(1\le h\le P_4-P_2=1\). Since \(\mathcal S_4=\mathcal A_4=\mathbb C\bar x\), we have \(\theta^2=\lambda\bar x\) with \(\lambda\ne0\). Thus \(x\) is a unit at \(T\), proving the nondivisibility of its powers. If \(v_\alpha(x/s^4)\in E_5\), then \(\Delta_4(x)\in sR_5\), whereas the Euler identity gives \(\overline{\Delta_4(x)}=-4\theta\bar x\ne0\), a contradiction. Finally, at a visible component of multiplicity \(k\ge2\), \(\theta\) has order \(k-1\); hence \(\theta^2=0\) in the generic local ring of \((s)\). The relation \(\theta^2=\lambda\bar x\) prevents \(x\) from being a unit there.
\(\square\)

\hypertarget{local-time-and-alternating-jacobians}{%
\hypertarget{section-3.3}{}\subsection{3.4. Local time and alternating Jacobians}\label{local-time-and-alternating-jacobians}}

We now describe the local form of the derivation near the canonical divisor. The functions under consideration are already meromorphic, so their convergent expansions are available before we apply any differential identity.

\needspace{3\baselineskip}

\hypertarget{result-3.6}{}\ResultHeading{Lemma 3.7 (local time).} At a general point of a transverse component
of multiplicity \(\ell\), there are analytic coordinates
\((\rho,y_1,y_2)\) and an analytic function \(\tau=\rho^{\ell+1}\) such that
\[
\alpha=dy_1\wedge dy_2,\qquad
s=d\tau\wedge dy_1\wedge dy_2,\qquad v=\partial_\tau.      \tag{3.6}\label{eq:3.7}
\]
A function with a positive pole of order \(n\) in \(\rho\) has derivative
of pole order \(n+\ell+1\). At a tangent component,
\[
\operatorname{ord}_T(vf)\ge\operatorname{ord}_T(f)-\ell.   \tag{3.7}\label{eq:3.8}
\]
If \(K_W\) is big, at least one component of the actual canonical divisor is transverse to \(\alpha\).

\ProofHeading{Proof.} The kernel of the closed rank-two form \(\alpha\) is
transverse to the divisor. The analytic flow theorem \cite[Theorem 1]{Carrillo2021}, together with closedness of \(\alpha\), gives coordinates
\((x,y_1,y_2)\) with \(\alpha=dy_1\wedge dy_2\) and
\(s=x^\ell h\,dx\wedge dy_1\wedge dy_2\), where \(h\) is a unit. Set
\[
\tau=\int_0^x z^\ell h(z,y_1,y_2)\,dz=x^{\ell+1}H.
\]
The unit \(H\) has an analytic \((\ell+1)\)-st root after shrinking.
Replacing \(x\) by \(\rho=xH^{1/(\ell+1)}\) gives \eqref{eq:3.7}.
The pole assertion follows by differentiation. In the tangent
case \(w\) is tangent, so \(\operatorname{ord}(w(f))\ge \operatorname{ord}(f)\); division by \(g\) gives \eqref{eq:3.8}.

To prove the last assertion, we use bigness to write \(K_W\equiv\varepsilon H+E\), with \(H\)
ample, \(\varepsilon>0\), and \(E\ge0\). The closed semipositive
\((2,2)\)-form \(\alpha\wedge\bar\alpha\) pairs positively with \(H\)
and nonnegatively with every effective divisor. Its pairing with
\(K_W\) is therefore positive. Decomposing \(\operatorname{div}_W(s)\)
shows that some restriction of \(\alpha\) is nonzero. \(\square\)

A meromorphic function has a convergent Laurent expansion in \(\rho\),
equivalently a convergent Puiseux expansion in \(\tau\), with analytic
transverse coefficients after shrinking. Its exponents belong to
\((\ell+1)^{-1}\mathbb Z\). The dominant-balance and resonance method goes back to Kowalevski \cite[Section 1]{Kowalevski1889}; see \cite[Section 3]{ConteMusette2003} for the coefficient-recursion formulation. A resonance is an exponent at which the
coefficient of the next unknown Laurent coefficient in a differential
equation vanishes. We use the term only for the recurrence of an
existing expansion. In particular, a first nonintegral exponent must
be tested against that recurrence.

\needspace{3\baselineskip}

\hypertarget{result-3.7}{}\ResultHeading{Lemma 3.8 (prime-degree field and Jacobian).} Suppose
\(K=\mathbb C(f_1,f_2,f_3)\subset L=\mathbb C(W)\),
\([L:K]\) is prime, and the \(f_i\) are algebraically independent.
Then
\[
J=\frac{df_1\wedge df_2\wedge df_3}{s}\notin K,\qquad
L=K(J).                                                 \tag{3.8}\label{eq:3.9}
\]
The three \(f_i\) cannot all have integral Laurent expansions in
the actual time \(\tau\) at a transverse component. If
\(v(f_1)=f_2,\ v(f_2)=f_3,\ v(f_3)=P(f_1,f_2,f_3)\), then
\[
\frac{v(J)}J=\frac{\partial P}{\partial f_3}.             \tag{3.9}\label{eq:3.10}
\]

\ProofHeading{Proof.} Algebraic independence gives \(J\ne0\). If \(J\in K\),
the form
\(\eta=(df_1\wedge df_2\wedge df_3)/J\) is a rational top form on
a projective space with function field \(K\), and its pullback is
\(s\). Resolve the rational map from \(W\) to that projective space.
The pullback of \(s\) remains holomorphic. A pole of order \(d\ge1\)
along a prime divisor downstairs would pull back to order
\(-ed+e-1<0\) at a prime divisor dominating it, where \(e\) is the
ramification index. Thus \(\eta\) has no divisorial pole. It would
be a nonzero holomorphic top form on \(\mathbb P^3\), a
contradiction. Hence \(J\notin K\), and the prime degree gives
\(L=K(J)\).

If the \(f_i\) have integral \(\tau\)-expansions, their derivatives
with respect to \(\tau,y_1,y_2\) do also. Formula \eqref{eq:3.7} expresses
\(J\) as their Jacobian determinant, so \(J\) also has integral
\(\tau\)-exponents. Every element of \(L=K(J)\) would then have
valuation divisible by \(\ell+1\), contradicting the normalized
divisorial value group \(\mathbb Z\). This includes \(\ell=1\).

More generally, if \(v(f_i)=F_i(f_1,f_2,f_3)\), then
\[
v(J)/J=\sum_{i=1}^3\partial F_i/\partial f_i.
\]
For the last assertion, we apply \(L_v\) to \(df_1\wedge df_2\wedge df_3=Js\).
Since \(L_vs=0\), only the third term in the derivative of the
wedge survives, giving \eqref{eq:3.10}. \(\square\)

\needspace{3\baselineskip}

\hypertarget{result-3.8}{}\ResultHeading{Corollary 3.9 (the quadratic involution).} Suppose the rational
map to \(\mathbb P^3\) defined by \(f_1,f_2,f_3\) has degree two.
Its involution negates every holomorphic two-form and three-form
on \(W\). Consequently it preserves the vector field \(v\), and
\(v(K)\subset K\).

\ProofHeading{Proof.} An invariant rational differential form descends through
a separable Galois field extension. A descended form cannot have
a divisorial pole on a smooth model of the quotient if its pullback
is holomorphic. In local parameters, the normal term of a pole of
order \(d\ge1\) pulls back with order \(-ed+e-1<0\), while tangential
terms have order at most \(-ed<0\); separability of the residue
extension prevents their cancellation. Resolve the rational map
to \(\mathbb P^3\) to apply this at every divisor. The descended
form would therefore be holomorphic on \(\mathbb P^3\), which has
no holomorphic forms of positive degree. Thus the invariant
parts vanish. Since both \(s\) and \(\alpha\) are negated,
\(\iota_v s=\alpha\) shows that \(v\) is invariant. \(\square\)

\needspace{3\baselineskip}

\hypertarget{result-3.9}{}\ResultHeading{Lemma 3.10 (local differential tests).} Use the actual time of \hyperlink{result-3.6}{Lemma 3.7}.

\begin{enumerate}
\def\labelenumi{(\alph{enumi})}
\item
  If \(v\) preserves a rational field \(K=\mathbb C(f_1,f_2,f_3)\) and
  \(v(f_i)=F_i(f_1,f_2,f_3)\), then its divergence
  \(\sum_i\partial F_i/\partial f_i\) has time-pole at most one. In particular,
  if this divergence is \(A(t)\in\mathbb C[t]\) and \(t\) has time-pole greater
  than one at a transverse component, then \(A\) is constant.
\item
  Consider an existing Puiseux expansion satisfying a differential
  equation with constant coefficients, or coefficients integral in time.
  Suppose its leading exponents are integers and its coefficient recurrence
  is triangular. If the indicial polynomial has no positive nonintegral
  rational root, the expansion has only integral exponents. The same holds
  for a linear system whose indicial matrix is invertible at every such
  shift. A first-order equation \(u'=a(\tau)u+b(\tau)\) with integral Laurent
  coefficients has this property if \(a\) has at most a simple pole with
  integral residue; a pole of \(a\) of order greater than one also excludes a
  first fractional exponent.
\item
  Suppose \(t,t',t''\) are algebraically independent and generate a
  rational subfield of prime index in \(\mathbb C(W)\). They cannot satisfy
  \[
  t'''=a\,tt'+b\,t^2+c\,t''+d\,t'+e\,t+f,
  \qquad a,b,c,d,e,f\in\mathbb C.
  \]
\end{enumerate}

\ProofHeading{Proof.} Part (a) follows from the logarithmic derivative identity in
\hyperlink{result-3.7}{Lemma 3.8}: for \(J=\tau^\lambda u\), with \(u\) a Puiseux unit,
\(J'/J=\lambda/\tau+u'/u\) has no pole of order greater than one.
For (b), take the first nonintegral shift. Every nonlinear or lower-order
contribution at that shift uses an earlier nonintegral coefficient, so
the recurrence reduces to its nonzero indicial coefficient times the
chosen coefficient. This is impossible. For the first-order assertion
the corresponding coefficient is \(q-\operatorname{res}a\); if \(a\) has a
higher pole, its product with the first fractional term has uniquely
lowest order.

For (c), take a transverse component, which exists by \hyperlink{result-3.6}{Lemma 3.7}. If \(t\)
is bounded, a first fractional term has an unmatched third derivative.
At a pole the possibilities and indicial polynomials are
\[
\begin{array}{c|c|c|c}
\text{condition}&\text{pole}&\text{leading coefficient}&I(q)\\\hline
a\ne0&2&12/a&(q+1)(q-4)(q-6)\\
a=0,\ b\ne0&3&-60/b&(q+1)(q^2-13q+60).
\end{array}
\]
If \(a=b=0\), a pole is impossible and the equation is linear. Part (b)
gives integral expansions in every case, contradicting \hyperlink{result-3.7}{Lemma 3.8}.
\(\square\)

We shall also use the same recurrence test for \(t''=kt^2+\) lower-order
terms at a double pole, with indicial polynomial \((q+1)(q-6)\), and for
\(t'''=k(t')^2+\) lower-order terms at a simple pole, with indicial
polynomial \((q+1)(q-1)(q-6)\). These polynomials are obtained by
substituting \(t=T\tau^{-p}+\varepsilon\tau^{q-p}\) in the displayed
dominant equation. Only actual meromorphic expansions are used.

\needspace{4\baselineskip}

\hypertarget{graded-subrings-and-odd-canonical-sections}{%
\hypertarget{section-3.4}{}\subsection{3.5. Graded subrings and odd canonical sections}\label{graded-subrings-and-odd-canonical-sections}}

For a finitely generated positively graded complex algebra \(B\) of
dimension \(n\), define
\[
e(B)=\lim_{q\to1^-}(1-q)^nH_B(q),\qquad
H_B(q)=\sum_{m\ge0}(\dim B_m)q^m.
\]
A polynomial ring of weights \(w_1,\ldots,w_n\) has multiplicity
\(1/(w_1\cdots w_n)\); a hypersurface of degree \(d\) has multiplicity
\(d/(w_1\cdots w_{n+1})\). Asymptotic Riemann--Roch gives \(e(R)=V\).

\needspace{3\baselineskip}

\hypertarget{result-3.10}{}\ResultHeading{Lemma 3.11 (field index and multiplication).} Let \(B\subset A\)
be finitely generated positively graded domains over \(\mathbb C\), of the
same dimension, and put
\(\delta=[\operatorname{Frac}A:\operatorname{Frac}B]\). Then
\(e(A)\ge\delta e(B)\). In particular, for a graded subring \(B\subset R\)
of dimension four,
\[
V\ge\delta e(B).                                        \tag{3.10}\label{eq:3.11}
\]
If \(U,V\) are nonzero finite-dimensional complex subspaces of a field
extension of \(\mathbb C\), then
\[
\dim(UV)\ge\dim U+\dim V-1.                              \tag{3.11}\label{eq:3.12}
\]

\ProofHeading{Proof.} We choose \(\delta\) homogeneous elements of \(A\)
that are independent over \(\operatorname{Frac}B\). Such a choice is
possible because homogeneous elements span the fraction-field extension.
Their \(B\)-multiples form a direct sum of shifted copies of \(B\) in
\(A\). Comparing Hilbert series as \(q\to1^-\), we obtain the first
inequality. Taking \(A=R\) gives \eqref{eq:3.11}. This argument does not require
\(A\) to be finite as a \(B\)-module.

For \eqref{eq:3.12}, we work in the finitely generated field generated by the two
subspaces and choose a valuation with one-dimensional leaves, for
example from a complete flag on a smooth projective model. We can choose
bases with distinct ordered values \(a_1<\cdots<a_e\) and
\(b_1<\cdots<b_f\). The products with values
\[
a_1+b_1<\cdots<a_e+b_1<a_e+b_2<\cdots<a_e+b_f
\]
are linearly independent, which proves \eqref{eq:3.12}. \(\square\)

\needspace{3\baselineskip}

\hypertarget{result-3.11}{}\ResultHeading{Lemma 3.12 (Cartier-index integrality).}
Let \(Y\) be a normal projective threefold with isolated singularities,
let \(D\) be an integral \(\mathbb Q\)-Cartier Weil divisor, and suppose
that \(rD\) is Cartier. For every integral Weil divisor \(T\), we have
\(rD^2T\in\mathbb Z\). If \(D\) is nef and \(T\) is effective, this
integer is nonnegative. In particular, for a terminal minimal threefold
with \(rK_Y\) Cartier, we have \(rV\in\mathbb Z\).

\ProofHeading{Proof.} We choose an ample Cartier divisor \(H\). For sufficiently
large \(k\), general members of \(|rD+kH|\) and \(|kH|\) are smooth and
avoid the finite singular set. The restrictions of \(D\) and \(T\) are
integral Cartier divisors on these surfaces. The difference of their
intersection numbers is \(rD^2T\), proving integrality. Nefness gives
nonnegativity when \(T\ge0\). For the final assertion, take \(D=T=K_Y\); effectivity is not needed for integrality. \(\square\)

\needspace{3\baselineskip}

\hypertarget{result-3.12}{}\ResultHeading{Lemma 3.13 (a prime canonical divisor).} Suppose \(p_g>0\) and \(rV=1\), where \(rK_Y\) is Cartier. Then the divisor of every nonzero canonical section on the canonical model is reduced and prime. The algebra \(A=R/(s)\) is a domain and, with \(a_m=P_m-P_{m-1}\),

\[
a_{i+j}\ge a_i+a_j-1\qquad(a_i,a_j>0).
\tag{3.12}\label{eq:3.13}
\]

\ProofHeading{Proof.} We write \(\operatorname{div}_Y(s)=\sum m_iT_i\). \hyperlink{result-3.11}{Lemma 3.12} gives \(rK_Y^2T_i\in\mathbb{Z}_{\ge0}\). A component surviving on the canonical model has positive intersection, because \(K_Z\) is ample. Thus \(1=rV=\sum m_i(rK_Y^2T_i)\) leaves precisely one surviving component, of multiplicity one. The prime-canonical-divisor criterion now gives the domain assertion. Explicitly, divisibility of a homogeneous section by \(s\) is equivalent to vanishing along that prime divisor: division extends as a reflexive canonical section on the normal variety \(Z\). This checks the actual graded pieces, without any normality assumption on the divisor or its affine cone. \hyperlink{result-3.10}{Lemma 3.11} proves \eqref{eq:3.13}. \(\square\)

Suppose \(\mathbb C(W)/K\) is quadratic with \(K\) rational. Let \(\iota\) be
its involution and set \(j|_{R_m}=(-1)^m\iota^*\). Then \(j(s)=s\), and a
section is \(j\)-invariant precisely when its normalized function belongs
to \(K\). We call a section odd if \(j\) acts on it as \(-1\).

We shall also use the reflexive canonical module of the canonical ring. It satisfies
\(\omega_R\simeq R(1)\). Its generator \(\Omega\) has degree \(-1\) and is odd:
on the smooth canonical cone,
\(\Omega=t_0^{-2}\eta\wedge dt_0\), where \(t_0\) is the fibre coordinate on
\(K_Z^{-1}\). The natural lift of \(\iota\) fixes this form, while the grading
involution negates it.

\needspace{3\baselineskip}

\hypertarget{result-3.13}{}\ResultHeading{Lemma 3.14 (canonical forms and odd sections).}

\begin{enumerate}
\def\labelenumi{(\alph{enumi})}
\item
  Let \(C\subset R^j\) be a finitely generated graded domain with the
  same fraction field. If a nonzero homogeneous rational top form \(\eta\)
  on \(\operatorname{Spec}C\) extends to a resolution of a proper birational
  model, then \(R\) has an odd section of degree \(\deg\eta+1\).
\item
  If \(C\) is Gorenstein with canonical singularities and its canonical
  generator has degree \(d\), then there is an odd section of degree \(d+1\).
\item
  Let \(B\subset R^j\) be polynomial, with the same fraction field and
  sum of generator weights \(S_B\). If a reduced homogeneous fraction
  \(G/F\in R^j\) has a homogeneous polynomial coordinate \(L\) of weight \(d\)
  as a factor of \(F\), then there is an odd section of degree \(S_B+1-d\).
\end{enumerate}

\ProofHeading{Proof.} For (a), at every codimension-one point of \(\operatorname{Spec}R\)
the map to \(\operatorname{Spec}C\) lifts to the proper resolution by the
valuative criterion. Pulling back \(\eta\) gives a regular reflexive form
\(h\Omega\). It is invariant, so \(h\) is odd, and its degree is \(\deg\eta+1\).
Part (b) follows because a canonical generator extends to a resolution.
Neither argument requires that \(R\) be Cohen--Macaulay or finite over \(C\).

For (c), multiply \(G/F\) by \(F/L\), divide \(G\) by the coordinate \(L\), and
subtract the quotient. The resulting subring has equation
\(LU=g(b_1,b_2,b_3)\ne0\). Such a hypersurface is normal and canonical:
successive general slices through \(L=U=0\) have Du Val equation
\(LU=t^n\), and inversion of adjunction \cite[Theorem 5.50(1)]{ref22} applies twice.
Its canonical generator has degree \(S_B-d\). Part (b) proves the claim.
\(\square\)

\needspace{4\baselineskip}

\hypertarget{two-forms-and-pencils}{%
\hypertarget{section-3.5}{}\subsection{3.6. Two-forms and pencils}\label{two-forms-and-pencils}}

\needspace{3\baselineskip}

\hypertarget{result-3.15}{}\ResultHeading{Proposition 3.15 (pencil threshold).} Suppose a rational
fibration has \(Q\ge_{\mathbb Q}pF\). For \(m\ge5\), if
\(p>4/(m-4)\), then \(\varphi_m\) is birational.

\ProofHeading{Proof.} We write \(Q\sim_{\mathbb Q}pF+E_p\), \(E_p\ge0\),
and put \(a_i=m-1-i/p\), \(i=1,2\).
For \(b_0=p/(1+p)\) the assumption gives
\(a_1b_0>3\), \(a_2b_0>2\). Choose rational
\(b<b_0\) close enough that these inequalities remain true.
\hyperlink{result-2.1}{Proposition 2.1} gives \(D\equiv bA+H_b\), \(H_b\ge0\).

For two distinct general fibres \(F_1,F_2\), put
\(L_2=(m-1)Q-2E_p/p\). Then
\(L_2-F_1-F_2\equiv a_2Q\) is nef and big. Vanishing gives
simultaneous restriction of \(K_W+\lceil L_2\rceil\).
On each fibre the target contains an adjoint with rational
part numerically \(a_2D\), and is nonzero by
\cite[Lemma 2.14]{ref8}, since \(a_2b>2\).
This subsystem lies in \(|mK_W|\), so the fibres are separated.

Similarly \(L_1=(m-1)Q-E_p/p\) restricts to an adjoint
\(|K_F+\lceil N\rceil|\), with \(N\equiv a_1D\).
Set \(t=a_1b>3\). Since \(D,A\) are nef and
\(D-bA\) effective, \(N^2\ge t^2A^2\ge t^2>9\).
A curve through a very general point outside the fixed support
of \(H_b\) has \(NC\ge tAC\ge t\).
But
\[
\frac4{1+\sqrt{1-8/t^2}}<t.
\]
Equation \eqref{eq:2.9} proves birationality on the general fibre.
Together with separation, it proves the proposition. \(\square\)

For clarity, the surface canonical-degree estimate used here and
in the subsequent reductions is the following published result:
if \(C\) passes through a very general point of a smooth surface
of general type and \((K_S^2,p_g(S))\ne(1,2)\), then
\(K_S C\ge2\) \cite[Lemma 2.5]{ref7}. Without that exception hypothesis,
nefness, bigness, and integrality give \(K_SC\ge1\).

\needspace{3\baselineskip}

\hypertarget{result-3.16}{}\ResultHeading{Lemma 3.16 (irregular fibres).} Suppose \(f:W\to B\) has general fibre
\(F\) with \(q(F)=1\), and some nonzero holomorphic two-form \(\alpha\) restricts
to zero on \(F\). With \(D=Q|_F\) and \(A=\sigma^*K_S\), one has \(AD\ge1\).
If \(f\) is the complete connected four-canonical pencil, then
\[
V\ge1/20,\qquad \varphi_{13}\text{ is birational}.
\]

\ProofHeading{Proof.} The form descends to a reflexive form on the terminal model.
Functorial pullback \cite[Theorems 1.3 and 5.2]{ref21} makes its restriction to
every exceptional divisor zero: the image has dimension at most one.
Locally \(\alpha=dt\wedge\beta\), where \(\beta|_F\) is a nonzero Albanese
differential. Hence all curves in \(E_\pi|_F\) are Albanese-vertical.

The actual canonical divisor \(K_F=D+E_\pi|_F\) has a nonzero horizontal
part \(H\), since it has positive intersection with a general Albanese
fibre. Every component of \(H\) survives on \(S\) and has positive integral
\(A\)-degree. Indeed a canonical-null curve has negative square by Hodge
index, whereas a horizontal curve has normalization of genus at least
one and hence nonnegative square by adjunction. Therefore \(D\ge H\) and
\(AD\ge AH\ge1\).

For the four-pencil, \hyperlink{result-2.1}{Proposition 2.1} gives \(D\ge_{\mathbb Q}A/5\) and
\(V\ge D^2/4\ge AD/20\). For birationality take \(N=8D\) in \eqref{eq:2.8}.
Then \(N^2\ge64/5\) and \(NC\ge16/5\) for a curve through a very general
point, by \cite[Lemma 2.5]{ref7}. Since
\[
\frac4{1+\sqrt{1-8/(64/5)}}<\frac{16}{5},
\]
Langer's criterion gives birationality on \(F\); the four-pencil separates
the outer fibres. \(\square\)

\needspace{4\baselineskip}

\hypertarget{small-plurigenera-and-an-eighth-degree-plateau}{%
\hypertarget{section-3.7}{}\subsection{3.7. Small plurigenera and an eighth-degree plateau}\label{small-plurigenera-and-an-eighth-degree-plateau}}

We next isolate a consequence of the saturation and Euler identities. It uses only the displayed plurigenera and the nonzero action of the chosen two-form on the pencil coordinate; no condition on the dimension of the full space of two-forms is needed.

\needspace{3\baselineskip}

\hypertarget{result-3.19}{}\ResultHeading{Proposition 3.17.}

Assume \(q(W)=0\) and

\[
P_1=P_2=P_3=1,\qquad P_4=P_5=2,\qquad P_6=P_7=4.
\]

Fix a nonzero holomorphic two-form \(\alpha\). Choose \(x_4\in R_4\) independent of \(s^4\), and assume \(v_\alpha(t)\ne0\) for \(t=x_4/s^4\). Then

\[
\theta^2=\lambda\bar x_4\ne0
 \qquad\text{for some }\lambda\in\mathbb{C}^*.
\tag{3.13}\label{eq:3.19}
\]

\ProofHeading{Proof.}
Corollary~\hyperlink{result-canonical-polar}{3.6} gives \(\theta^2\ne0\).
By \hyperlink{result-3.4}{Proposition 3.3},
\(\mathcal S_4=\mathcal A_4=\mathbb C\bar x_4\), so
\(\theta^2=\lambda\bar x_4\) with \(\lambda\ne0\). \(\square\)

\needspace{3\baselineskip}

\hypertarget{result-3.20}{}\ResultHeading{Corollary 3.18.}

Under the hypotheses of \hyperlink{result-3.19}{Proposition 3.17}, put \(k=24\lambda\). Then

\[
r_6=\Delta_4(x_4)\in R_6,\qquad
 z_6=\frac{\Delta_6(r_6)-kx_4^2}{s^2}\in R_6.
\tag{3.14}\label{eq:3.21}
\]

There is a component \(T\) of multiplicity one with \(\alpha|_T\ne0\). On every such component \(x_4\) is a unit at the generic point, whereas on every other component of \(\operatorname{div}_Z(s)\) its order is at least the multiplicity of \(s\).

\ProofHeading{Proof.}

The Euler relation gives \(\overline{\Delta_6(r_6)}=24\theta^2\bar x_4=k\bar x_4^2\). The difference is in \(sR_7=s^2R_6\), proving \eqref{eq:3.21}. If \(\alpha|_T=0\), then \(\theta=0\) in the full generic Artin ring. If \(\alpha|_T\ne0\) and the multiplicity is \(a\geq2\), its square vanishes there because \(2(a-1)\geq a\). Equation~\eqref{eq:3.19} gives the asserted vanishing of \(x_4\) on all these components. At a component of multiplicity one with \(\alpha|_T\ne0\), the element \(\theta\) is a unit, so \(x_4\) is a unit as well. If there were no such component, \(x_4/s\) would have no divisorial pole and would belong to \(R_3=\mathbb{C} s^3\), contrary to the choice of \(x_4\). \(\square\)

\needspace{3\baselineskip}

\hypertarget{result-3.21}{}\ResultHeading{Lemma 3.19 (an eighth-degree plateau).} Assume \(q(W)=0\) and
\[ (P_1,\ldots,P_7)=(1,1,1,2,2,4,4). \]
Fix a nonzero holomorphic two-form \(\alpha\), and choose \(t=x/s^4\) with \(v_\alpha(t)\ne0\). If \(L_6\) is rational and \([\mathbb C(W):L_6]\) is prime, then \(P_8=P_9=6\) is impossible.

\ProofHeading{Proof.} \hyperlink{result-3.19}{Proposition 3.17} and \hyperlink{result-3.20}{Corollary 3.18} give
\(\theta^2=\lambda\bar x\ne0\). Put \(k=24\lambda\),
\[
r=t',\qquad z=t''-kt^2\in E_6.
\]
There is a reduced transverse component where \(x\) is a unit.
Its actual time is \(\tau=\rho^2\), and
\[
t=(6/k)\tau^{-2}+\cdots,\qquad t^2\notin E_6.              \tag{3.15}\label{eq:3.22}
\]
The second Euler identity and \(R_9=sR_8\) give
\[
f''-2ktf\in E_8\qquad(f\in E_6).                          \tag{3.16}\label{eq:3.23}
\]
If \(1,t,r,z\) form \(E_6\), either
\(z'\in E_6+\mathbb Ct^2\), giving
\[
t'''=2ktt'+A t''+B t'+Ct+Dt^2+E,
\]
or \(z'\) completes that space to \(E_8\), and \eqref{eq:3.23} gives
\[
t''''-2k(t')^2-4ktt''+2k^2t^3
 =A+Bt+Ct'+D(t''-kt^2)+Et^2+F(t'''-2ktt').
\]
At the leading pole \eqref{eq:3.22}, the indicial polynomials are respectively
\[
(q+1)(q-4)(q-6),\qquad
(q+1)(q-1)(q-6)(q-8).
\]
A first nonintegral exponent has no forcing from earlier integral
exponents and cannot have a nonzero coefficient. Thus the three
sixth coordinates have integral actual-time expansions.

If \(z\in\langle1,t,r\rangle\), then
\(t''=kt^2+A+Bt+Ct'\), whose indicial roots are \(-1,6\).
Choose \(u\) completing \(E_6\). Either its derivative lies in
\(E_6+\mathbb Ct^2\), giving a first-order linear equation for \(u\),
or it completes \(E_8\), and \eqref{eq:3.23} gives
\[
u''=2ktu+A+Bt+Ct'+Du+Et^2+Fu'.
\]
The principal exponents in the latter equation satisfy
\(p(p-1)=12\), namely \(-3,4\). In either case, the already integral
forcing and integral homogeneous exponents exclude a first
fractional exponent of \(u\). Prime-degree Jacobian reconstruction
contradicts the resulting integral expansions of all three
generators. This proves the lemma. \(\square\)

\needspace{4\baselineskip}

\hypertarget{reduction-to-a-marked-canonical-pencil}{%
\needspace{10\baselineskip}\hypertarget{section-4}{}\section{4. Reduction to a marked canonical pencil}\label{reduction-to-a-marked-canonical-pencil}}

We retain the notation of Sections 2 and 3. In particular, \(Y\) is a terminal minimal model, \(Z\) is its canonical model, \(Q\) is the pullback of the actual divisor of a canonical section \(s\), and \(E_m=R_m/s^m\). A prime component of \(\operatorname{div}_Z(s)\) is called transverse to \(\alpha\) when the tangential restriction of \(\alpha\) is nonzero. All divisor-order arguments on the canonical model are made at its smooth codimension-one points.

\needspace{4\baselineskip}

\hypertarget{strictness-and-the-square-two-endpoint}{%
\hypertarget{section-4.1}{}\subsection{4.1. Strictness and the square-two endpoint}\label{strictness-and-the-square-two-endpoint}}

\needspace{3\baselineskip}

\hypertarget{result-4.1}{}\ResultHeading{Lemma 4.1 (strictness at the curve threshold).}

Let \(D\) be nef and big on a smooth surface, and let \(C\) be a general member of a free moving curve system. Suppose \(D\equiv\beta C+E\), \(E\ge0\), \(\beta>0\), and a linear system contains \(|K_F+\lceil N_a\rceil|\), for an actual rational divisor \(N_a\equiv aD\). Put \(t=a-\beta^{-1}>0\) and \(\xi=DC\). Its restriction to \(C\) is birational if either \(t\xi>2\), or \(t\xi\ge2\) and \(D^2>\beta\xi\).

\ProofHeading{Proof.}

The Zariski decomposition of \(tD+E/\beta\) has positive part \(P\equiv tD+B\), where \(B,N\ge0\), \(B+N=E/\beta\), and \(PN=0\). This follows from maximality of the positive part, since \(tD\) is nef. Thus \(PC\ge t\xi\). In the second case \(DE>0\). If \(DB=0\), Hodge index gives \(B^2<0\) for nonzero \(B\), whereas \(PB=B^2\ge0\). If \(B=0\), the equality \(PN=0\) instead contradicts \(DE>0\). Hence \(DB>0\). If \(BC=0\), Hodge index gives \(B^2\le0\), while

\[
B(B+N)=B(D/\beta-C)=DB/\beta>0.
\]

This forces \(BN>0\), contradicting \(PN=tDN+BN=0\). Therefore \(PC>2\).

The negative part depends only on the numerical class. Subtract it from the actual divisor and put \(P'=N_a-C-N\). This is nef and big, numerically equivalent to \(P\), and \(\lceil P'\rceil+C\le\lceil N_a\rceil\). Its fractional support is independent of the general member \(C\). Resolve that support by \(\rho:F'\to F\) and choose \(C\) avoiding the finitely many centres; \(\rho\) is then an isomorphism near \(C\). Vanishing gives \(H^1(F',K_{F'}+\lceil\rho^*P'\rceil)=0\), so the adjoint system with \(C\) added restricts onto \(|K_C+T|\), with \(\deg T\ge\lceil PC\rceil\ge3\). Sections descend by their codimension-one divisor inequalities, and their restrictions near \(C\) are unchanged. Thus the original system contains this complete restriction, which has degree at least \(2g(C)+1\) and is very ample. \(\square\)

\needspace{3\baselineskip}

\hypertarget{result-4.2}{}\ResultHeading{Proposition 4.2.}

Suppose a four-canonical pencil has connected rational base and its general fibre has a canonical genus-two pencil. If \(m\ge11\), \(\xi\ge2/(m-10)\), and \(\varphi_m\) is not birational, then

\[
\xi=\frac2{m-10},\qquad D^2=\frac\xi5,\qquad V=\frac\xi{20}.
\]

\ProofHeading{Proof.}

The canonical restriction lemma gives \(D\ge_{\mathbb{Q}}A/5\ge_{\mathbb{Q}}C/5\). Apply the adjoint restriction lemma and \hyperlink{result-4.1}{Lemma 4.1} with \(a=m-5\), \(\beta=1/5\). Separation of the surface and curve pencils then forces the first two equalities. On the special model,

\[
\xi\le DA\le DK_F\le5D^2=\xi.
\]

\hyperlink{result-2.2}{Proposition 2.2}, with \((a,b)=(4,1)\), gives the volume. \(\square\)

\needspace{3\baselineskip}

\hypertarget{result-4.3}{}\ResultHeading{Proposition 4.3 (an impossible surface endpoint).}

On the special pencil model of \hyperlink{result-2.2}{Proposition 2.2}, the conditions

\[
4K_Y\sim F_Y+Z_Y,\qquad Z_Y\ge0,\qquad
 A^2=2,\quad D^2=1/10,\quad DA=1/2
\]

cannot hold simultaneously. The divisor relation is an integral Weil-divisor relation.

\ProofHeading{Proof.}

Put \(H=G|_F+\mu^*Z_Y|_F\). Then \(5D\sim_{\mathbb{Q}}K_F+H\) and \(D(K_F-A+H)=0\). Every exceptional curve of \(\sigma\) occurs in \(K_F-A\), so \(D\) descends to a nef big rational divisor \(D_0\) on \(F_0\). If \(T=\sigma_*H\),

\[
5D_0\sim_{\mathbb{Q}}K_{F_0}+T,\qquad D_0T=0,\qquad K_{F_0}T=1/2.
\]

The support of \(T\) is negative definite. For a component \(C\) with \(k=K_{F_0}C>0\) and coefficient \(t_C\), adjunction gives

\[
0=(K_{F_0}+T)C\ge k+t_CC^2\ge k-t_C(k+2),
 \qquad t_Ck\ge\frac{k^2}{k+2}.
\tag{4.1}\label{eq:4.1}
\]

There is exactly one such component \(C_0\), with canonical degree one, coefficient \(1/2\), and square \(-3\); adjunction makes it smooth rational. All other components are smooth rational \((-2)\)-curves. If \(B\) is their positive Cartan matrix and \(v\) the incidence vector with \(C_0\), their coefficients satisfy \(Bb=v/2\) and \(v^tB^{-1}v=1\). The positive definite \((-2)\)-blocks are ADE. For \(A_n\) the inverse diagonal is \(j(n+1-j)/(n+1)\), for \(D_n\) it is \(1,2,\ldots,n-2,n/4,n/4\), and its minimum for \(E_6,E_7,E_8\) is \(4/3,3/2,2\), respectively. Two nonzero entries in one connected block give norm larger than one. The remaining possibilities are two \(A_1\) blocks, one \(A_3\) block attached at its middle, or one \(D_n\) block attached at a vector vertex. Their coefficients belong to \(\{1/4,1/2\}\). These configurations are simple normal crossing trees: larger intersection multiplicities or a multiple crossing would violate their intersection matrices.

The divisor \(K_F+H-\sigma^*(K_{F_0}+T)\) is exceptional and numerically trivial, hence zero. Factor \(\sigma\) into point blowups. Intermediate crepant boundaries remain effective, as pushforwards of \(H\). A new coefficient is the sum at the centre minus one. The only possible effective blowup of the listed boundary is a crossing of two \(1/2\)-components, with new coefficient zero. Consequently all positive coefficients of \(H\) are \(1/4\) or \(1/2\), and \((F,H)\) is klt.

We next show \(Z_Y=0\). Restrictions of distinct horizontal prime divisors to a very general fibre are reduced and share no component. A strict component of \(Z_Y\) would contribute a positive integer to \(H\), which is impossible. The same is true of exceptional divisors with curve centre on \(Y\): \(Y\) is smooth at that generic centre, so their discrepancies and pullback orders have the corresponding integrality. All positive components of \(H\) therefore map to points of \(Y\). The effective restriction \(\mu^*Z_Y|_F\) has degree zero on each of its components by projection formula. Negative definiteness makes it zero.

Take a general very ample smooth surface on \(Y\). On it the effective divisors induced by \(F_Y\) and \(Z_Y\) have zero mutual intersection, and their positive weighted sum is a nef big multiple of \(K_Y\). Their squares are nonnegative; Hodge index forces one to be numerically zero. The divisor induced by \(F_Y\) is nonzero, so the one induced by \(Z_Y\) must vanish. Hence \(Z_Y=0\). This also rules out fixed divisors in special fibres, which restriction to a general fibre alone would not detect.

We explain the Cartier descent needed to finish. Since \((F,G|_F)\) is klt, the proof of \cite[Lemma 3.4(1)]{ref13} applies to \eqref{eq:2.4}: inversion of adjunction gives plt near \(F\), connectedness gives plt near \(F_Y\), and terminality deals with its complement. Thus \(F_Y\) is normal and klt. As \(Y\) has only isolated singularities, there is no divisorial different, and

\[
K_F+G|_F=\mu_F^*K_{F_Y}.
\]

\(4(K_F+G|_F)\) is integral Cartier with degree zero on every exceptional curve. Over the local ring of a point of \(F_Y\), rational surface singularities give \(H^1(\mathcal O)=0\). \cite[Theorem (12.1)(i)]{ref24} makes this line bundle trivial there. Pushing down a local rational trivialization proves that \(4K_{F_Y}\) is Cartier.

For a plt pair, restriction of a \(\mathbb{Q}\)-Cartier Weil divisor to the integral boundary preserves its local Cartier index \cite[Appendix A]{ref18}. This is the actual Weil restriction class here, since \(Y\) is smooth in codimension one along \(F_Y\). Adjunction and \(F_Y\sim4K_Y\) give \(K_{F_Y}\sim5K_Y|_{F_Y}\), so \(20K_Y\) is Cartier along \(F_Y\); off it, \(4K_Y\) is Cartier. Thus \(20K_Y\) is globally Cartier, whereas \hyperlink{result-2.2}{Proposition 2.2} gives \(V=1/40\), contradicting \hyperlink{result-3.11}{Lemma 3.12}. \(\square\)

\needspace{4\baselineskip}

\hypertarget{the-geometric-reduction}{%
\hypertarget{section-4.2}{}\subsection{4.2. The geometric reduction}\label{the-geometric-reduction}}

\needspace{3\baselineskip}

\hypertarget{result-4.4}{}\ResultHeading{Proposition 4.4.} If \(p_g(X)>0\) and \(\varphi_{14}\) is not birational, then
\[p_g=P_2=P_3=\chi(\mathcal O_X)=1,\qquad q=0,\qquad P_4=2.\]
The four-canonical pencil has connected base \(\mathbb P^1\). Its general minimal fibre has type \((K^2,p_g,q)=(1,2,0),(2,2,0)\) or \((2,1,1)\). In the regular cases the canonical moving pencil has genus two and every nonzero global two-form restricts nontrivially. In the irregular case every global two-form restricts to zero. In all cases \(h^{2,0}(X)=1\).

\ProofHeading{Proof.} Hu--Yan's theorem \cite{ref19} excludes \(P_3\ge2\), the irregular fifth-canonical theorem \cite{ref11} excludes \(q>0\), and Chen--Zuo's theorem \cite{ref12} excludes \(\chi(\mathcal O_X)\le0\). Hence \(p_g=P_2=P_3=1\), \(q=0\), and \(h^{2,0}=\chi>0\). Proposition~\hyperlink{result-polar-inequality}{3.5} gives \(P_4-P_2\ge h\), and hence \(P_4\ge2\).

Suppose that the four-canonical image has dimension at least two. Let \(F\) be a general member of its free moving part \(M_4\), and \(C\) a general irreducible member of the restricted system. Then the restriction parameters satisfy \(\mu,\beta\ge1/4\). Writing \(A=\sigma^*K_{F_0}\) and \(D=Q|_F\), adjunction and maximality of the positive part of the surface Zariski decomposition give
\[A\ge D+M_4|_F\ge\tfrac54M_4|_F.\]
Thus \(A^2\ge2\), and the moving-curve inequality \cite[Lemma 2.4]{ref7} gives \(AC\ge2\). Canonical restriction yields \(D\ge_{\mathbb Q}A/5\), so \(\xi=DC\ge2/5\). The separation lemmas \cite[Lemmas 3.1--3.2]{ref14} apply in degree fourteen. If \(g(C)\ge3\), the numerical curve criterion gives \(9\xi\ge2g(C)-2\ge4\), and \((13-4-4)\xi>2\). A three-dimensional image has \(C^2>0\), hence \(g(C)\ge3\). The only remaining case has image dimension two, \(g(C)=2\), and \(\xi=2/5\).

Here \(D^2\ge\xi/4=1/10\). Strict inequality gives birationality by the strict curve criterion above, with the degree-fourteen adjoint of numerical part \(9D\) and \(\beta=1/4\). At equality,
\[\tfrac12=D^2+\xi\le DA\le5D^2=\tfrac12.\]
Hodge index gives \(A^2=2\). A general subpencil of \(|M_4|\) has an integral relation \(4K_Y\sim F_Y+Z_Y\), \(Z_Y\ge0\). This is the impossible surface endpoint proved above. Thus the four-canonical image is a curve.

Its Stein base is rational since \(q=0\). If \(P_4\ge3\), its moving part contains at least two connected fibres, giving \(Q\ge_{\mathbb Q}F/2\) and birationality in degree thirteen by \hyperlink{result-3.15}{Proposition 3.15}. Thus \(P_4=2\) and the complete moving pencil consists of connected fibres. The same polar inequality now gives \(0<h\le P_4-P_2=1\), so \(h=\chi=1\).

If \(A^2\ge3\), the degree-fourteen adjoint has \(N=9D\), with
\[N^2\ge\tfrac{81}{25}A^2>9,\qquad NB\ge\tfrac95AB\ge\tfrac{18}5\]
for curves \(B\) through a very general point. Langer's bound is less than three, a contradiction. Hence \(A^2\le2\). For a regular fibre the Hodge restriction lemma gives \(p_g(F)\ge h^{2,0}+1\ge2\); Noether's inequality leaves \((1,2),(2,2),(2,3)\). For an irregular fibre, Debarre's inequality \cite{ref16} \(K_{F_0}^2\ge2p_g(F)\), the positive restriction of \(s\), and \(\chi(\mathcal O_F)>0\) leave \((2,1,1)\).

For \(A^2=2\), a canonical moving curve of genus at least four gives \(10\xi\ge2g(C)-2\ge6\), which suffices in degree fourteen. If its genus is three, the rounded curve inequalities give
\[m\xi\ge4+\lceil(m-10)\xi\rceil.\]
\hyperlink{result-2.4}{Lemma 2.4} gives \(\xi\ge1/2\). Strictness again handles \(\xi>1/2\), and at equality the strict curve criterion reduces to \(D^2=1/10\), \(DA=1/2\), excluded by the surface endpoint lemma. A \((2,3)\) surface has a free canonical system with general curve of genus three: Bombieri's canonical-system theorem \cite[Theorem 9]{ref3} excludes a pencil and a fixed part; its canonical image is a plane and its degree is at least two, so \(K^2=2\) leaves no base point. This removes \((2,3)\).

The assertions about restrictions of two-forms follow from \hyperlink{result-3.2}{Lemma 3.2}(c). \(\square\)

A minimal \((1,2)\) surface has no fixed canonical curve and one simple canonical base point. Indeed, writing \(K=M+Z\) gives \(KM=1\), \(KZ=0\) and \(M^2\le1\). If \(M^2=0\), the pencil is free and adjunction would give the odd number \(2g-2=1\). Thus \(M^2=1\) and there is one simple base point; then \(MZ=KZ=0\) and \(Z^2=0\), so Hodge index gives \(Z=0\). Blowing up that point gives a genus-two pencil. For a regular \((2,2)\) surface whose canonical moving pencil has genus two, the pencil is free. If \(m_i\) were its base-point multiplicities, then \(2=KM+\sum m_i\) and \(M^2=\sum m_i^2\), contradicting \(2M^2\le(KM)^2\) when a base point exists. Thus \(K=C+Z\), with \(C^2=0\), \(KC=2\), \(KZ=0\) and \(Z^2=-2\).

\needspace{4\baselineskip}

\hypertarget{fibres-with-canonical-square-two}{%
\hypertarget{section-4.3}{}\subsection{4.3. Fibres with canonical square two}\label{fibres-with-canonical-square-two}}

The irregular \((2,1,1)\) case is excluded by \hyperlink{result-3.16}{Lemma 3.16}. It remains
to treat a regular \((2,2,0)\) fibre whose canonical pencil has genus two.

For a minimal regular \((2,2,0)\) fibre with canonical genus-two pencil, write
\[
K_{F_0}\sim C+Z,\qquad
C^2=0,\quad K_{F_0}C=2,\quad K_{F_0}Z=0,\quad Z^2=-2.
\tag{4.2}\label{eq:4.2}
\]
The moving pencil is free. Indeed, if its moving part \(M\) had base-point multiplicities \(m_i\), adjunction on the resolution would give \(2=K_{F_0}M+\sum m_i\) and \(M^2=\sum m_i^2\). A base point would imply \(K_{F_0}M\le1\) and \(M^2\ge1\), contrary to Hodge index \(2M^2\le(K_{F_0}M)^2\).

The fixed cycle \(Z\) is a single fundamental ADE cycle. Each component is a canonical \((-2)\)-curve. A further \((-2)\)-curve meeting its support would have negative intersection with the nef moving divisor \(C\). Thus the support is a union of complete ADE configurations. Moreover \(Z E_i=-C E_i\le0\). In a finite ADE root lattice, an effective anti-nef integral cycle of square \(-2\) is the fundamental root cycle; \(Z^2=-2\) permits only one connected block. Its coefficients are at most six, with coefficients at least five occurring only in type \(E_8\). We use the intersection configurations of rational double points as listed in \cite[Chapter III]{ref1}. This observation is not needed for the irregular exclusion.

\needspace{3\baselineskip}

\hypertarget{result-4.5}{}\ResultHeading{Proposition 4.5.} Let \(S\) be a smooth minimal surface of general type with \((K_S^2,p_g(S),q(S))=(2,2,0)\), whose canonical moving system is a genus-two pencil. Write \(K_S\sim C+Z\) with \(Z\) fixed. Then \(Z\) is the reduced fundamental cycle of an \(A_{2r-1}\) configuration for some integer \(r\ge1\).

\ProofHeading{Proof.} The discussion above shows that \(Z\) is one ADE fundamental cycle. Let \(\iota\) be the fibrewise hyperelliptic involution and \(g:S\to\mathbb P^1\) the canonical pencil. It extends to an automorphism of the minimal surface. Use its natural action on forms: \(H^0(K_S)\) is negative and \(H^0(2K_S)\) is positive.

Set \(V_n=g_*\omega_{S/\mathbb P^1}^{\otimes n}\). The degree formulas \cite[Remark 2.11]{ref5} give \(\deg V_1=4\) and \(\deg V_2=14\). The injection \(\operatorname{Sym}^2V_1\to V_2\) has torsion cokernel of length \(14-12=2\), denoted \(\mathcal O_\tau\) in \cite{ref5}. \cite[Proposition 4.8]{ref5} gives
\[
V_3^+\simeq\det V_1\otimes\mathcal O_{\mathbb P^1}(\tau)
\simeq\mathcal O_{\mathbb P^1}(6).
\]
Twisting by \(\omega_{\mathbb P^1}^{\otimes3}\) shows \(h^0(3K_S)^+=1\). Riemann--Roch and vanishing give \(P_2(S)=5\), \(P_3(S)=9\), and thus \(h^0(3K_S)^-=8\).

A fundamental cycle of a rational double point satisfies \(H^1(Z,\mathcal O_Z)=0\). Adjunction and duality give
\[
H^0(Z,\mathcal O_Z(K_S+Z))=H^0(Z,\omega_Z)=0.
\]
Hence \(h^0(K_S+Z)=2\). The base-point-free pencil trick for \(|C|\) shows that the kernel of
\[
H^0(K_S)\otimes H^0(2K_S)\longrightarrow H^0(3K_S)
\]
has dimension \(h^0(2K_S-C)=h^0(K_S+Z)=2\). The image has dimension eight, is negative, and vanishes along \(Z\). It equals the entire negative tricanonical eigenspace.

Let \(\rho:S\to\bar S\) be the canonical contraction and \(p=\rho(Z)\). Bombieri's theorem \cite[Theorem 2]{ref3} gives global generation of \(3K_S\). Therefore a positive tricanonical section is nonzero at \(p\). Equivariant evaluation gives positive character on \(\omega_{\bar S,p}^{\otimes3}\), hence on \(\omega_{\bar S,p}\). Averaging a local generator produces an invariant nowhere-zero local generator \(e\) of \(\omega_{\bar S}\).

The two canonical sections vanish simultaneously only at \(p\). Write them locally as \(ue\) and \(ve\). They are negative, so \(u,v\) are odd. Their local intersection length is
\[
\operatorname{length}\mathcal O_{\bar S,p}/(u,v)=K_{\bar S}^2=2.
\]
A rational double point has embedding dimension three. The length-two quotient has embedding dimension one, so \(u,v\) have independent cotangent classes. Choose a third eigen-coordinate \(w\). The length-two quotient also ensures a nonzero \(w^2\) coefficient in the local hypersurface equation. Equivariant Weierstrass preparation and completion of the square give
\[
w^2=f(u,v).
\tag{4.3}\label{eq:4.3}
\]
The relation has positive character. The residue formula \(du\wedge dv/(2w)\), together with the positive canonical character, forces \(w\) to be even. Consequently \(f(-u,-v)=f(u,v)\).

A rational-double-point branch germ in double-plane form has multiplicity at most three, by the ADE normal forms \cite[Chapter III]{ref1}. Thus the even germ \(f\) has multiplicity two. A rank-two quadratic part gives \(A_1\). In rank one, the equivariant analytic Morse lemma gives \(f=u'^2+h(v')\) with \(u',v'\) odd and \(h\) even. Isolation makes the order of \(h\) finite, say \(2r\). This is an \(A_{2r-1}\) singularity, whose fundamental cycle is reduced. \(\square\)

\needspace{3\baselineskip}

\hypertarget{result-4.6}{}\ResultHeading{Lemma 4.6.} Under the hypotheses of \hyperlink{result-4.5}{Proposition 4.5}, the fixed cycle \(Z\) has type \(A_1,A_3\) or \(A_5\).

\ProofHeading{Proof.} We write its type as \(A_{2r-1}\). For \(r\ge2\), the identity \(C E_i=-Z E_i\) equals one on the two endpoints and zero on all interior components. Thus the endpoints are sections of the canonical pencil, whereas the connected interior chain \(A_{2r-3}\) is vertical. No other canonical \((-2)\)-curve meets \(Z\), by nefness of \(C\). The relative canonical contraction therefore contracts exactly this interior block to a point \(q\) of type \(A_{2r-3}\).

The positive tricanonical section \(z\) from \hyperlink{result-4.5}{Proposition 4.5} is nonzero throughout \(Z\). Near \(q\), the involution acts positively on the canonical generator, whereas relative degree-one canonical sections have negative character. All these sections vanish at \(q\). It is therefore one of the distinguished points of the genus-two relative canonical algebra \cite[Lemmas 4.1 and 4.3, Theorem 4.7]{ref5}.

Let \(\ell\) be the coefficient of the Horikawa divisor \(\tau\) at this fibre. Since \(\deg\tau=2\), one has \(\ell\le2\). The local algebra in \cite[proof of Lemma 4.1]{ref5} has degree-two relation
\[
t^\ell y=Q_2(x_0,x_1,t),
\]
where \(Q_2\) is quadratic in \(x_0,x_1\) and its specialization at \(t=0\) has rank at least one. At \(q\) use the chart \(y=1\). The degree-three generator \(z\) is a unit, so the degree-six equation eliminates it. Locally the surface equation is exactly
\[
t^\ell-Q_2(x_0,x_1,t)=0.
\tag{4.4}\label{eq:4.4}
\]

For \(\ell=1\) this is smooth. For \(\ell=2\), a rank-two quadratic part in \(x_0,x_1\) gives \(A_1\). In rank one, a \(t\)-dependent invertible linear change puts it in the form
\[
t^2-x_0^2-c(t)x_1^2=0,\qquad c(0)=0.
\]
If \(\operatorname{ord}_t c\ge2\), the equation factors analytically, since
\(1-c(t)x_1^2/t^2\) is a unit with a square root. This contradicts normality. Therefore \(\operatorname{ord}_t c=1\). The analytic Morse lemma in \(x_0,t\) leaves a nonzero quartic term in \(x_1\), giving \(A_3\). Comparing with \(A_{2r-3}\) yields \(r\le3\). The case \(r=1\) was already allowed. \(\square\)

We use multiplier ideals on surfaces with rational double points in the following explicit form. For a log resolution \(\rho:\widetilde S\to S\) and an effective Cartier divisor \(B\),
\[\mathcal J(S,cB)=\rho_*\mathcal O_{\widetilde S}(K_{\widetilde S}-\rho^*K_S-\lfloor c\rho^*B\rfloor).\]
A local function \(h\) belongs to this ideal exactly when
\[\operatorname{ord}_E(h)+A_S(E)>c\operatorname{ord}_E(B)\tag{4.5}\label{eq:4.5}
\]
for every divisorial valuation, where \(A_S(E)\) is the log discrepancy. If a local double cover has ramification function \(r\), the canonical-divisor formula gives
\[A_T(v)+v(r)=eA_U(v_0)\tag{4.6}\label{eq:4.6}
\]
for a valuation upstairs with normalized restriction \(v_0\) and ramification index \(e\). These are equalities on models carrying the divisors and are preserved by crepant resolutions.

\needspace{3\baselineskip}

\hypertarget{result-4.7}{}\ResultHeading{Proposition 4.7.} Let \(S\) be a smooth minimal surface of general type with \((K_S^2,p_g(S),q(S))=(2,2,0)\) and a canonical genus-two pencil. For every effective divisor \(A_0\in|K_S|\),
\[
|3K_S\otimes\mathcal J(S,A_0/5)|
\quad\text{is birational}.
\tag{4.7}\label{eq:4.7}
\]

\ProofHeading{Proof.} Let \(z\) span the positive tricanonical eigenspace of \hyperlink{result-4.5}{Proposition 4.5}. We first prove that \(z\) belongs to the indicated multiplier ideal. The verification is local, using a local canonical generator to view a tricanonical section as a function.

At the canonical point \(p\) obtained by contracting \(Z\), equation \eqref{eq:4.3} describes a double plane branched along two smooth branches \(B_1,B_2\) with contact order \(r\le3\). A selected canonical divisor is the pull-back of a smooth plane germ \(L\): its equation is a nonzero linear combination of the original parameters \(u,v\). The pair
\[
\left(\mathbb C^2,\tfrac12 B_1+\tfrac12 B_2+\tfrac15 L\right)
\tag{4.8}\label{eq:4.8}
\]
is klt. Here is the local check. Blow up the common tangent cluster. While all three strict transforms pass through the same centre, at most \(r\) successive steps are possible and the exceptional boundary coefficient is at most \(r/5\le3/5\). After \(B_1\) and \(B_2\) separate, at most one coefficient-\(1/2\) branch continues with \(L\). The sum \(1/2+1/5<1\) makes further exceptional coefficients decrease. All strict-transform coefficients are below one. The same calculation covers coincidence of \(L\) with one branch, whose coefficient is then \(7/10\). Thus every boundary coefficient on a resolution is less than one. The double-cover formula makes \((\bar S,A_0/5)\) klt near \(p\). The canonical resolution is crepant, so the multiplier ideal is trivial near \(Z\).

Away from \(Z\), use the relative canonical model \(T\) and its double cover \(\psi:T\to U\) of the conic bundle in \cite[Theorem 4.7]{ref5}. Choose a nonvanishing global canonical section \(a\). All distinguished points of this double cover lie over \(Z\): at such a point both relative degree-one generators vanish, so all global canonical sections vanish. They can therefore be omitted on this open set.

The remaining conic charts are described in \cite[Remark 4.14]{ref5}. They are smooth with a reduced fibre or twice a smooth fibre component, or an additional \(A_1\) point on a double fibre. The latter pair has local equation
\[
x^2=t y,\qquad \text{fibre }t=0.
\]
Its quasi-étale cover is the smooth plane with \(t=u^2\), \(y=v^2\), \(x=uv\). Thus in every case the pair \((U,\tfrac15 U_t)\) is klt; its fibre threshold is at least \(1/2\).

The ratio \(r_0=z/a^3\) is a local ramification generator for \(\psi\). There is no missing zero from the outer base: \(V_3^+\otimes\omega_{\mathbb P^1}^{\otimes3}\simeq\mathcal O_{\mathbb P^1}\), and its global generator \(z\) is nowhere zero as a generator of this line bundle. On a chart where \(a\) is nonzero the relative canonical algebra expresses the cover as \(r_0^2=b\). Locally \(A_0\) is the pull-back of a fibre divisor, or is absent. For a divisorial valuation over this chart, \eqref{eq:4.6} gives
\[
A_T(v)+v(r_0)
=eA_U(v_0)>\tfrac15 e\,v_0(U_t)=\tfrac15v(A_0).
\]
The valuation criterion \eqref{eq:4.5} and crepancy of the resolution give \(z\in\mathcal J(S,A_0/5)\) away from \(Z\) as well.

Finally, \(s_{A_0}H^0(2K_S)\) belongs to the same ideal automatically. Its map is the bicanonical map, which has degree exactly two here: it recovers the pencil base, and the restriction to a general canonical curve is the complete bicanonical system, because
\[
h^0(2K_S-C)=h^0(K_S+Z)=2,\qquad P_2(S)=5.
\]
The target on the curve has dimension three. The character of \(s_{A_0}H^0(2K_S)\) is negative, whereas that of \(z\) is positive. Their ratios therefore supply a function outside the invariant quadratic subfield. Adding \(z\) makes the map birational. \(\square\)

For a log resolution \(\sigma:F\to S\) of \(A_0\), the defining formula for the multiplier ideal identifies the sections in \eqref{eq:4.7} with
\[
H^0\left(F,K_F+\left\lceil\tfrac95\sigma^*A_0\right\rceil\right),
\tag{4.9}\label{eq:4.9}
\]
using \(\lceil(9/5)\sigma^*A_0\rceil=2\sigma^*A_0-\lfloor\sigma^*A_0/5\rfloor\). For an arbitrary smooth birational model, membership in the multiplier ideal still gives inclusion into \eqref{eq:4.9}, by the divisorial inequalities. Equality of pushforwards is not asserted without a log resolution.

\needspace{3\baselineskip}

\hypertarget{result-4.8}{}\ResultHeading{Proposition 4.8.} Under the assumptions of \hyperlink{result-2.2}{Proposition 2.2}, suppose that the general minimal fibre has invariants \((K^2,p_g,q)=(2,2,0)\) and its canonical pencil has genus two. Then \(\varphi_{14}\) is birational.

\ProofHeading{Proof.} On a common sufficiently fine smooth model, the adjoint restriction formula and the actual inequality \(5D\ge\sigma^*A_0\) of \hyperlink{result-2.2}{Proposition 2.2} give
\[
|14K_W||_F
\succeq |K_F+\lceil9D\rceil|
\succeq
\left|K_F+\left\lceil\tfrac95\sigma^*A_0\right\rceil\right|.
\]
The last system is birational by \hyperlink{result-4.7}{Proposition 4.7} and \eqref{eq:4.9}. The four-canonical pencil is contained in \(|14K_W|\) and separates the outer fibres. Birationality on the general fibre therefore proves birationality on \(W\). \(\square\)

\needspace{4\baselineskip}

\hypertarget{marked-curves-and-the-low-plurigenera}{%
\hypertarget{section-4.4}{}\subsection{4.4. Marked curves and the low plurigenera}\label{marked-curves-and-the-low-plurigenera}}

\hyperlink{result-4.4}{Proposition 4.4} and the preceding fibre exclusions reduce failure of birationality in degree fourteen to a rational four-canonical pencil of \((1,2)\) surfaces. We establish the lower bound \(\xi\ge1/3\) on its marked canonical curve. Throughout this subsection,
\[p_g=P_2=P_3=\chi=1,\qquad q=0,\qquad P_4=2.\]
The common restriction results give \(D\ge_{\mathbb Q}A/5\) and \(V\ge D^2/4\ge\xi/20\).

We use the normalized baskets and Reid's Riemann--Roch formula. The global terminal index is the least common multiple of the basket denominators, including the deformation baskets of nonquotient terminal points; see \cite[Section 2.1, Remark 2.1]{ref6}.

For each basket pair set

\[
\gamma=\min\{r-2b,b\},\quad d=b-\gamma,\quad
 e=\frac b2-\frac{b^2}{r}-\frac\gamma6.
\]

Substitution of the low-plurigenus assumptions in the Reid Riemann--Roch formula gives

\[
\sum b=14,\qquad\sum\gamma=5,\qquad\sum d=9,
 \qquad\sum e=\frac16-V.
\tag{4.10}\label{eq:4.10}
\]

Indeed the equations for \(P_2,P_3\) give \(V=8-\sum b+\sum b^2/r\) and \(\sum b=14\); the correction at \(3b=r\) in the \(P_4\) equation is \((3b-r)_+=b-\gamma\). Every \(e\) is nonnegative. For \(r\ge3b\) it equals \(b/3-b^2/r\); for \(2b<r<3b\) it equals \(\gamma d/[6(3\gamma+2d)]\).

\needspace{3\baselineskip}

\hypertarget{result-4.9}{}\ResultHeading{Lemma 4.9 (index bound).} If \(p_g=P_2=P_3=\chi=1\) and \(P_4=2\), every basket denominator is at most \(25\).

\ProofHeading{Proof.} The pair \((1,2)\) is harmless. If \(r\ge3b\), then \(b=\gamma\le5\) and \(e<1/6\) give
\(r<6b^2/(2b-1)\le50/3\). Otherwise \(r=3\gamma+2d\), with
\(1\le\gamma\le5\), \(d\le9\), and \(d(\gamma-2)<3\gamma\).
For \(\gamma=1,2,3,4,5\) the maximum denominators are respectively
\(21,24,25,22,23\). \(\square\)

\needspace{3\baselineskip}

\hypertarget{result-4.10}{}\ResultHeading{Lemma 4.10.} In the residual \((1,2)\) four-pencil,
\(\xi\ge5/17\) and \(V\ge1/68\).

\ProofHeading{Proof.} Apply \hyperlink{result-2.5}{Lemma 2.5} to the published seventeenth-canonical
theorem \cite[Theorem 1.3(i)]{ref14}. Then \(V\ge\xi/20\ge1/68\). \(\square\)

We record the numerical identities before improving \(\xi\). They do not require a birational pluricanonical map.

\needspace{3\baselineskip}

\hypertarget{result-4.11}{}\ResultHeading{Lemma 4.11.}

Under the low-plurigenus assumptions, put \(u=P_5\), \(v=P_6\). Then \(u\in\{2,3\}\),

\[
P_7=u+v-2,\qquad
 B^{(5)}=(8-2u+v)(1,2)+(1+2u-v)(2,5)
 +(1-u+v)(1,3)+(3-u)(1,4),
\tag{4.11}\label{eq:4.11}
\]

and \(V\le(u+2v-7)/60\). When \(u=2\), one has \(3\le v\le5\).

\ProofHeading{Proof.}

Use the general initial-basket formulas \cite[(3.6), (3.8), (3.10)]{ref9}, before their subsequent specialization \(P_2=0\). In our case they read \(n^0_{1,2}=9\), \(n^0_{1,3}=u+2\) and \(n^0_{1,4}=3-u-\sigma_5\), where \(\sigma_5\ge0\) counts the initial entries of denominator at least five. Since \(u\ge2\), we have \(u=2\) or \(3\). For \(u=3\), \(\sigma_5=0\). For \(u=2\), \cite[(3.10)]{ref9} reads \(0=v-P_7-\epsilon\) with \(\epsilon\ge0\); monotonicity \(P_7\ge v\) gives \(\epsilon=\sigma_5=0\). In either case \cite[(3.10)]{ref9} gives the first identity in \eqref{eq:4.11}; \cite[(3.8)]{ref9} gives the basket. Its volume is \(\sum b^2/r-6=(u+2v-7)/60\). Every nontrivial prime packing decreases volume by the prime-packing formula. For \(u=2\), positivity of volume gives \(v\ge3\), and nonnegativity of the \((2,5)\) coefficient gives \(v\le5\). \(\square\)

Following \cite[Definition 3.3]{ref14}, define

\[
\begin{aligned}
 u_m&=\dim\operatorname{Im}\bigl(H^0(W,M_m)\longrightarrow
                         H^0(F,M_m|_F)\bigr),\\
 v_{m,-j}&=\dim\operatorname{Im}\bigl(H^0(F,M_m|_F-jC)\longrightarrow
                         H^0(C,M_m|_C)\bigr)\quad(j\ge0),
 \qquad v_m=v_{m,0}.
\end{aligned}
\]

Here \(H^0(W,M_m)\) is identified with \(R_m\) after removing the fixed part, whereas the source defining \(v_{m,-j}\) is the full space of sections on \(F\). If \(H_m=h^0(F,M_m|_F)\), the restriction sequences give

\[
u_m\le H_m,\qquad h^0(F,M_m|_F-C)=H_m-v_m,\qquad
 h^0(F,M_m|_F-2C)=H_m-v_m-v_{m,-1}.
\]

We shall use \cite[Propositions 3.6--3.7]{ref14} with the separation parameter \(m_2=5\): the fourth pencil distinguishes outer fibres, and the adjoint restriction lemma shows that degree six distinguishes the canonical curves on \(F\). Here are the exact consequences needed. If \(M_a\ge F+S_1\) with \(S_1\) moving outside the outer pencil, and \(S_1|_F\) defines a different curve pencil, put \(d_1=S_1\cdot C\ge2\). \cite[Proposition 3.6(1.2)]{ref14} gives birationality in degree at most

\[
\left\lfloor\frac1\beta+\frac{2a}{d_1}
       +\frac1\mu\left(1-\frac2{d_1}\right)\right\rfloor+2
 \le\left\lfloor9+\frac{2a-8}{d_1}\right\rfloor+2.
\tag{4.12}\label{eq:4.12}
\]

If it is the same curve pencil, \cite[Proposition 3.6(2)]{ref14} gives

\[
(n+a+1)\xi\ge\lceil n\xi\rceil+2\quad(n\xi>1),
 \qquad \varphi_{n+a+1}\text{ birational if }n\xi>2.
\tag{4.13}\label{eq:4.13}
\]

The auxiliary condition \(n\ge m_2-a\) is automatic in every application below. If \(M_a|_F\ge C+C_1\) with a different moving curve family, \cite[Proposition 3.7]{ref14} gives birationality in degree at most \(a+6\) for \(a=7,8\): if \(C_1C=2\), use its part (i); if \(C_1C>2\), part (ii) gives the same bound.

\needspace{3\baselineskip}

\hypertarget{result-4.12}{}\ResultHeading{Proposition 4.12.}

In the residual \((1,2)\) pencil one has
\[
\xi\ge\frac13,\qquad V\ge\frac1{60}.
\tag{4.14}\label{eq:4.14}
\]

\ProofHeading{Proof.} Suppose \(\xi<1/3\). \hyperlink{result-4.10}{Lemma 4.10} gives \(\xi\ge5/17\) and
excludes birationality in every degree at most fifteen. If
\(D\ge_{\mathbb Q}C/3\), the rounded inequality in degree twelve gives
\(12\xi\ge2+\lceil4\xi\rceil\ge4\), a contradiction.

A two-fibre subsystem in degree at most seven is also impossible.
It gives \(Q\ge_{\mathbb Q}2F/7\), \(D\ge_{\mathbb Q}2A/9\), and
\(m\xi\ge2+\lceil(m-9)\xi\rceil\) for large \(m\). \hyperlink{result-2.4}{Lemma 2.4}
gives \(\xi\ge1/3\), a contradiction.

For \(P_5=3\), use the geometric alternatives in \cite[proof of Proposition
4.6, Case 1]{ref14}. A pencil in degree five contains two outer fibres. A
different curve pencil on \(F\) gives \(5\xi\ge2\), while two canonical
members give \(D\ge_{\mathbb Q}2C/5\). Otherwise the fifth restriction is
one canonical curve \(C_5\). On a general fifth surface \(S_5\), let
\(\widehat C\) be a generic irreducible member of its restricted system.
If \(C_5\widehat C>0\), this intersection is at least two and again
\(\xi\ge2/5\). If it is zero and the restriction contains \(2\widehat C\),
the changed flag has the same \(\xi\), with \(\mu\ge1/5\), \(\beta\ge2/5\).
Now \(\alpha(12)\ge(7/2)\xi\ge35/34>1\), so \hyperlink{result-2.3}{Proposition 2.3} gives
\(12\xi\ge4\). In the remaining irreducible restriction, the cited
Subcase 1.2 gives
\(\xi\ge(\sigma_5^*K_{S_{5,0}}\cdot\widehat C)/6\ge1/3\):
\(p_g(S_5)\ge2\) and \(S_5\) is not a \((1,2)\) surface. All cases contradict
the assumption.

Let \(P_5=2\), so \(3\le P_6=P_7\le5\) by \hyperlink{result-4.11}{Lemma 4.11}. The restriction
sequences in \cite[proof of Proposition 4.6, Case 2]{ref14} give the following
alternatives. If \(v_6\ge3\), then \(6\xi\ge4\). If \(v_{6,-1}\ge2\),
vanishing gives
\[
|12K_W||_F\succeq|K_F+\lceil D\rceil+C+C_1|,\qquad C_1C\ge2,
\]
which is birational after restriction to \(C\); both pencils are
separated. If \(h^0(F,M_6|_F)\ge4\), or this dimension is three and
\(v_6=1\), then \(M_6|_F\ge2C\), already excluded. Hence that dimension
is at most three.

If \(P_6\ge5\), or \(P_6=4\) and \(u_6\le2\), the kernel supplies
\(M_6\ge F+F_1\). The same outer pencil gives two fibres. A different
inner pencil gives birationality by degree thirteen through \eqref{eq:4.12}.
For the same inner pencil, \eqref{eq:4.13} at \(n=4\) gives
\(11\xi\ge2+\lceil4\xi\rceil\ge4\). If \(P_6=4\), \(u_6=3\), the
remaining restriction has dimension three and \(v_6=2\), so \(6\xi\ge2\).

Finally, \(P_6=3\) gives the ancestor
\(7(1,2)+2(2,5)+2(1,3)+(1,4)\) of volume \(1/60\). Its possible first
prime packings have volumes \(1/420,-1/60,1/120,1/210\), all below the
bound \(V\ge1/68\) from \hyperlink{result-4.10}{Lemma 4.10}. Thus there is no packing. Its Cartier
index is \(60\), and Reid's formula gives \(a_8=a_{16}=2\), where
\(a_m=P_m-P_{m-1}\). \hyperlink{result-3.12}{Lemma 3.13} contradicts \(a_{16}\ge2a_8-1\).
This proves \(\xi\ge1/3\), and \(V\ge\xi/20\) proves the volume estimate.
\(\square\)

\needspace{4\baselineskip}

\hypertarget{the-sixth-map-and-canonical-subrings}{%
\needspace{10\baselineskip}\hypertarget{section-5}{}\section{5. The sixth map and canonical subrings}\label{the-sixth-map-and-canonical-subrings}}

From now on, we suppose that \(\varphi_{14}\) is not birational. The geometric reduction leaves a four-canonical pencil of \((1,2)\) surfaces over \(\mathbb P^1\). Under this assumption, \hyperlink{result-4.12}{Proposition 4.12} gives \(\xi\ge1/3\) and \(V\ge1/60\). We use these bounds to control the sixth-canonical system and the possible baskets.

The curve restriction criterion applies with \(\mu=1/4\), \(\beta=1/5\) and
\[\alpha(m)=(m-10)\xi.\tag{5.1}\label{eq:5.1}\]
Since \(\xi\ge1/3\), its degree-fourteen restriction contains a complete genus-two line bundle of degree at least four. Such a bundle has a map of degree at most two; a bundle of degree at least five is very ample. Both pencil bases are recovered by the adjoint restriction lemma. Thus a nonbirational \(\varphi_{14}\) has degree exactly two. Its involution fixes both bases and acts hyperelliptically on the general canonical curve. Adjunction gives the sign \(-1\) on the relative canonical form and hence on \(s\) and \(\alpha\).

\needspace{4\baselineskip}

\hypertarget{two-fibre-endpoints-and-the-sixth-canonical-system}{%
\hypertarget{section-5.1}{}\subsection{5.1. Two-fibre endpoints and the sixth-canonical system}\label{two-fibre-endpoints-and-the-sixth-canonical-system}}

We first consider the numerical endpoint \(V=1/36\). Assume \(p_g=P_2=P_3=\chi=1\), \(P_4=2\), and the actual marked-curve condition with \(\xi=1/3\). We can enumerate the three possible baskets directly, as follows.

The gap sum is \(\sum e=5/36\), and \(36\mid r_X\). Thus some denominator is divisible by nine and some by four. Since all denominators are at most 25, the first is 9 or 18. Index 18 would have to be \((7,18)\), of \(\gamma=4,e=1/9\); the remaining gamma budget is one and remaining gap is \(1/36\), too small for any denominator divisible by four (with gamma at most one, such a pair has gap at least \(1/12\)). Therefore a pair \((4,9)\) occurs, with \(\gamma=1,d=3,e=1/18\).

There cannot be two such pairs: their remaining gap would be \(1/36\), again smaller than the minimum positive gap \(1/24\) of a pair with denominator divisible by four. After one \((4,9)\), the remaining gap is \(1/12\). A pair with denominator divisible by four and gap at most \(1/12\) must be one of
\[
 (1,4),\qquad(3,8),\qquad(5,12),
\]
with gaps \(1/12,1/24,1/12\). The only positive-gap pairs of gap at most \(1/24\) are \((2,5)\), of gap \(1/30\), and \((3,8)\), of gap \(1/24\). Therefore choosing \((3,8)\) forces a second \((3,8)\). Filling the zero-gap entries yields exactly
\[
 \begin{aligned}
 B_1&=6(1,2)+3(1,3)+(1,4)+(4,9),\\
 B_2&=3(1,2)+2(1,3)+(4,9)+(5,12),\\
 B_3&=4(1,2)+2(3,8)+(4,9).
 \end{aligned}
\]
The first two have \(r_XV=1\), and violate the domain inequality respectively by
\[
 (a_4,a_6,a_{10})=(1,2,1),\qquad
 (a_4,a_9,a_{13})=(1,2,1).
\]
For \(B_3\), the marked local index must be nine. By the terminal-type table \cite[Remark 2.1]{ref6}, the unique index-nine constituent comes from a \(cA/9\) point of axial weight one, whose index-one cover is smooth. \hyperlink{result-2.6}{Lemma 2.6} would give a linear canonical coefficient of character \(-4=5\) modulo nine, but the linear coordinate characters are \(1,8,4\). This is impossible. This proves the endpoint exclusion.

\needspace{3\baselineskip}

\hypertarget{result-5.1}{}\ResultHeading{Proposition 5.1.} In the \((1,2)\) row, suppose a subsystem of \(|6K_W|\) has moving part containing \(2F\). If \(\varphi_{14}\) fails to be birational, then
\[
\xi=\frac13,\qquad D^2=\frac1{12},\qquad V=\frac1{36}.
\]

\ProofHeading{Proof.} We write \(6Q\sim_{\mathbb Q}2F+E\), \(E\ge0\). Exact restriction gives \(D\ge_{\mathbb Q}A/4\ge_{\mathbb Q}C/4\). Vanishing gives in degree fourteen an actual surface adjoint with nef numerical part \(10D\), since \(13-6/2=10\). The curve threshold is therefore \((10-4)\xi=6\xi\ge2\). The outer and inner pencils are separated by the four-pencil restriction lemma. Strictness forces \(\xi=1/3\) and \(D^2=\xi/4=1/12\).

On the special pencil model, the integral relation \(6K_Y\sim2F_Y+Z\) gives
\[
4D\sim_{\mathbb Q}K_F+G|_F+\tfrac12\pi^*Z|_F,
\qquad G\ge0.
\]
Hence \(\xi\le DA\le4D^2=\xi\), and the special-model orthogonality lemma yields \(V=(2/6)D^2=1/36\). \(\square\)

The endpoint calculation above gives precisely three baskets and excludes each one. Consequently, within the regular \((1,2)\) row,
\[
6Q\ge_{\mathbb Q}2F\quad\Longrightarrow\quad\varphi_{14}\text{ is birational}.
\]

\needspace{3\baselineskip}

\hypertarget{result-5.2}{}\ResultHeading{Proposition 5.2.} Suppose \(D\ge_{\mathbb Q}C/3\) and \(\varphi_{14}\) is not birational. Then
\[
\boxed{\xi=\frac13,\qquad D^2=\frac19.}
\]

\ProofHeading{Proof.} Degree fourteen contains an adjoint with numerical part \(9D\). The strictness lemma with \(\beta=1/3\) gives exactly the two asserted equalities. \(\square\)

The equality here does not imply numerical proportionality of \(D\) and \(A\).

\needspace{3\baselineskip}

\hypertarget{result-5.3}{}\ResultHeading{Lemma 5.3.} Let \(\sigma:F\to S\) be the smooth minimal-model contraction on the special fibre. If it is not an isomorphism over the canonical base point \(p_S\), then \(\xi\ge2/5\). If \(\xi<2/5\), then
\[
A\cdot D=\xi.
\tag{5.2}\label{eq:5.2}
\]

\ProofHeading{Proof.} We factor \(\sigma\) into point blowups. If a blowup occurs over \(p_S\), consider the strict transform \(e_0\) of the exceptional curve of its first blowup. Every canonical member is smooth at the simple base point. Thus the actual canonical divisor on \(F\) has coefficient two along \(e_0\): one from the selected canonical member and one from the first discrepancy. Subsequent blowups do not change this coefficient on the strict transform.

A general resolved canonical member meets \(e_0\) transversely once, away from all later blowup centres, and misses every other component of the selected canonical divisor. Consequently the coefficient of \(D\) along \(e_0\) is \(\xi\). The actual inequality \(5D\ge K_F\) gives \(5\xi\ge2\).

If \(\xi<2/5\), there is therefore no blowup over \(p_S\). A general canonical member avoids all the other, finitely many blowup centres. Its strict transform on \(F\) is numerically \(A\). Resolving its remaining base point leaves the intersection with the pullback of \(D\) unchanged. This proves \eqref{eq:5.2}. \(\square\)

\needspace{3\baselineskip}

\hypertarget{result-5.4}{}\ResultHeading{Lemma 5.4.} Let \(S\) be a smooth minimal surface with \((K_S^2,p_g(S),q(S))=(1,2,0)\). If an integral curve \(B\) through a very general point satisfies \(K_S\cdot B=1\), then \(B\in|K_S|\).

\ProofHeading{Proof.} Negative curves are rigid and belong to at most countably many Hilbert components, so \(B^2\ge0\). Hodge index gives \(B^2\le1\), and adjunction makes \(B^2\) odd. Thus \(B^2=1\) and \(B\equiv K_S\).

If \(B\not\sim K_S\), the nontrivial numerically trivial line bundle \(\mathcal O_S(K_S-B)\) has no section. Riemann--Roch and \(\chi(\mathcal O_S)=3\) imply \(h^0(B)\ge3\). If \(|B|\) is composed with a pencil, its moving part contains at least twice a covering member, contradicting \(K_S\cdot B=1\). Otherwise, resolving the moving system gives a surface image with degree times generic degree at most \(B^2=1\). The image is a plane and the map is birational, contradicting general type. Hence \(B\sim K_S\). \(\square\)

On a common resolution, choose a very general point outside the support of an effective representative of \(D-A/5\). By this lemma and \hyperlink{result-2.1}{Proposition 2.1}, every integral curve \(B\) through that point satisfies
\[
9D\cdot B\ge
\begin{cases}
9\xi\ge3,& B\text{ is a canonical-pencil member},\\
18/5,& A\cdot B\ge2.
\end{cases}
\tag{5.3}\label{eq:5.3}
\]

Suppose first that \(D^2>1/9\). For \(N=9D\), one has \(N^2>9\), so the right side of the surface adjoint bound is strictly less than three. This contradicts \eqref{eq:5.3}. By the adjoint restriction formula the adjoint is contained in the fourteen-canonical restriction and the outer fibres are separated. Therefore failure of \(\varphi_{14}\) gives
\[
D^2\le\frac19.
\tag{5.4}\label{eq:5.4}
\]

There is exactly one component \(\Gamma\) of the actual canonical divisor on the special fibre having positive \(A\)-degree: its coefficient is one and \(A\cdot\Gamma=1\). Indeed \(A\cdot K_F=K_S^2=1\), while \(A\) is nef and Cartier. Write
\[
c=\operatorname{coeff}_\Gamma(D)=A\cdot D.
\tag{5.5}\label{eq:5.5}
\]
The equality follows from \(0\le D\le K_F\) and the vanishing of the \(A\)-degree of all the other components. On a resolution of the canonical pencil, \(A-C\) is effective, whence \(c\ge\xi\). Also \hyperlink{result-2.1}{Proposition 2.1} and nefness give
\[
D^2\ge\frac c5.
\tag{5.6}\label{eq:5.6}
\]
Thus \(1/3\le c\le5/9<1\).

The horizontal prime divisor containing \(\Gamma\) consequently has discrepancy \(1-c\in(0,1)\). It is exceptional over \(Y\). A divisor centred along a curve on a terminal threefold has positive integral discrepancy, since the generic point of that curve is smooth. A divisor over a Gorenstein terminal point likewise has integral discrepancy. Its centre is therefore a non-Gorenstein terminal point. The projection formula gives
\[
D\cdot\Gamma=0.
\tag{5.7}\label{eq:5.7}
\]

Assume now \(D^2=1/9\). Langer's bound is exactly three. By \eqref{eq:5.3} its obstruction curve must be a general canonical member and
\[
\xi=\frac13.
\]
\hyperlink{result-5.3}{Lemma 5.3} therefore gives \(A\cdot D=1/3\) on the special fibre. Since \(A^2=1\) and \(D^2=1/9\), equality in the Hodge index theorem yields
\[
D\equiv\frac13 A.
\]
But \eqref{eq:5.7} says \(D\cdot\Gamma=0\), whereas this numerical equivalence gives \(D\cdot\Gamma=1/3\). This contradiction proves
\[
\boxed{D^2<\frac19}.
\tag{5.8}\label{eq:5.8}
\]
Combining \eqref{eq:5.6}, \(c\ge\xi\ge1/3\), and \eqref{eq:5.8} proves \(1/15\le D^2<1/9\) and \(1/3\le\xi\le A\cdot D<5/9\).

If \(\xi<2/5\), \hyperlink{result-5.3}{Lemma 5.3} gives \(c=\xi\). The curve \(\Gamma\) is \(D\)-null, so \(\Gamma^2<0\) by Hodge index. Put \(T=5D-K_F\ge0\). Its coefficient along \(\Gamma\) is \(5c-1\). Intersections of distinct effective curves are nonnegative, and therefore
\[
-K_F\cdot\Gamma=T\cdot\Gamma\ge(5c-1)\Gamma^2.
\]
Adjunction gives
\[
2-2p_a(\Gamma)\ge(5c-2)\Gamma^2>0.
\]
Thus \(p_a(\Gamma)=0\), and \(\Gamma\) is smooth rational. This is a statement about the selected canonical member; the general canonical member still has genus two.

If the special model does blow up \(p_S\), \hyperlink{result-5.3}{Lemma 5.3} instead gives \(\xi\ge2/5\). Every curve in \eqref{eq:5.3} then has \(9D\)-intersection at least \(18/5\). For \(D^2>1/10\), the right side of the surface adjoint bound is less than \(18/5\), again a contradiction. Hence \(D^2\le1/10\) in this alternative. These arguments complete the geometric reduction on the special fibre.

On a sufficiently high smooth model, put \(R_m=H^0(W,mK_W)\) and
\[
\nu_6=\dim\operatorname{im}\bigl(R_6\longrightarrow H^0(F,6K_W|_F)\bigr).
\]

\needspace{3\baselineskip}

\hypertarget{result-5.5}{}\ResultHeading{Lemma 5.5.} Under the present counterexample hypothesis, the displayed restriction map has one-dimensional kernel.

\ProofHeading{Proof.} Its kernel is nonzero because the four-pencil, multiplied by \(s^2\), has a member containing the selected general fibre. If the kernel had dimension at least two, remove its fixed part and resolve to obtain a moving system \(|S_1|\) with
\[
M_6\ge F+S_1.
\]
If this is the same outer pencil, \(M_6\ge2F\), contrary to the two-fibre exclusion above.

Otherwise apply \cite[Proposition 3.6]{ref14} with
\[
m_0=4,\quad m_1=6,\quad m_2=5,\quad\mu=\frac14,\quad\beta=\frac15.
\]
The separating degree \(m_2+1=6\) is valid by the adjoint restriction formula: its restriction contains the canonical pencil and it separates outer fibres.

If the restricted moving system and the canonical pencil differ, their intersection degree \(\delta\) on a general genus-two member is at least two. Part (1.2) of that proposition gives birationality of \(\varphi_{n+1}\) for
\[
n\ge\max\left\{5,\left\lfloor5+\frac{12}{\delta}+4\left(1-\frac2\delta\right)\right\rfloor+1\right\}.
\]
The expression inside the floor is \(9+4/\delta\le11\), so \(n=12\) gives birationality in degree thirteen.

If the inner pencils agree, equation \eqref{eq:4.13}, successively with \(n=4\) and \(n=3\), gives
\[
11\xi\ge\lceil4\xi\rceil+2\ge4,
\qquad
10\xi\ge\lceil3\xi\rceil+2\ge4.
\]
The second application is allowed because the first gives \(3\xi>1\). Thus \(\xi\ge2/5\). Equation \eqref{eq:4.13} with \(n=6\) now gives birationality in degree \(6+6+1=13\), since \(6\xi>2\). Both conclusions contradict nonbirationality in degree fourteen. \(\square\)

\needspace{3\baselineskip}

\hypertarget{result-5.6}{}\ResultHeading{Proposition 5.6.} One has \(\nu_6\le3\) and \(P_6\le4\). If \(P_6=4\), then \([\mathbb C(W):L_6]=2\) and \(L_6=L_{14}\).

\ProofHeading{Proof.} We resolve the restricted system and let \(N_6\) be its free moving divisor. It contains the connected canonical-pencil base field, by the adjoint restriction formula.

If its image is a curve, it factors through that canonical pencil: the function field of its rational base is algebraically closed in the surface function field. Write \(N_6\sim aC\). A base-point-free subspace of dimension \(\nu_6\) has \(a\ge\nu_6-1\). If \(\nu_6\ge4\), then \(a\ge3\), so \(D\ge_{\mathbb Q}C/2\). The curve criterion gives
\[
(14-5-2)\xi=7\xi\ge\frac73>2,
\]
a contradiction.

If its image is a surface, its generic degree is at least two. Indeed all sixth-canonical ratios belong to the degree-fourteen field, and that field has index two in \(\mathbb C(W)\). The outer base field already belongs to both fields, so this remains a nontrivial degree-two obstruction on the generic outer fibre. A nondegenerate surface in \(\mathbb P^{\nu_6-1}\) has degree at least \(\nu_6-2\). Consequently
\[
36D^2\ge N_6^2\ge2(\nu_6-2).
\tag{5.9}\label{eq:5.9}
\]
The first inequality follows from \(6D-N_6\ge0\) and the nefness of both divisors. The strict inequality \(D^2<1/9\) rules out \(\nu_6\ge4\). Finally \hyperlink{result-5.5}{Lemma 5.5} gives \(P_6=\nu_6+1\le4\).

Suppose now that \(P_6=4\), so \(\nu_6=3\). A curve image would give \(N_6\sim aC\) with \(a\ge2\), hence \(D\ge_{\mathbb Q}C/3\). The endpoint refinement proved in \hyperlink{result-5.2}{Proposition 5.2} states that this would force
\[
\xi=\frac13,\qquad D^2=\frac19
\]
under nonbirationality of \(\varphi_{14}\). It contradicts \eqref{eq:5.8}. Thus the restricted image is a surface.

By \(\xi\le A\cdot D<5/9\), \(6\xi<10/3<4\). Formula \eqref{eq:2.10} places the sixth-canonical ratios on a general \(C\) in \(H^0(C,\mathcal O_C(3p))=\langle1,u\rangle\), where \(u\) is a hyperelliptic coordinate. Since the image is a surface and the inner base field is already present, this restriction is nonconstant and generates the hyperelliptic fixed field. The map has degree two on the generic \(C\), hence degree two on the generic \(F\) and on \(W\).

The inclusion of fields defined by \(\varphi_6\) and \(\varphi_{14}\), both of index two, is therefore an equality. This proves the common-field assertion without any assumption on \(P_7\). \(\square\)

\hypertarget{the-two-low-plurigenus-configurations}{%
\hypertarget{section-5.3}{}\subsection{5.2. The two low-plurigenus configurations}\label{the-two-low-plurigenus-configurations}}

\hypertarget{result-low-configurations}{}\ResultHeading{Proposition 5.7 (the remaining numerical configurations).}
Under fourteen-canonical failure,
\[
(P_5,P_6,P_7)=(2,4,4)\quad\text{or}\quad(3,4,5).
\]

\ProofHeading{Proof.}
Put \(u=P_5\) and \(v=P_6\). \hyperlink{result-4.11}{Lemma 4.11}
and \hyperlink{result-5.6}{Proposition 5.6} give
\(u\in\{2,3\}\), \(3\le v\le4\), and \(P_7=u+v-2\).
Corollary~\hyperlink{result-canonical-polar}{3.6} gives \(v>u\), excluding
\((u,v)=(3,3)\).
For \((u,v)=(2,3)\), Lemma~4.11 gives
\[
B^{(5)}=7(1,2)+2(1,3)+(1,4)+2(2,5),\qquad V(B^{(5)})=1/60.
\]
Its possible first nontrivial prime packings have volume losses
\(1/30,1/70,1/84,1/120\), all larger than
\(1/60-1/68=1/510\). The lower bound \(V\ge1/68\) from
\hyperlink{result-4.10}{Lemma 4.10} therefore permits no packing.
Thus \(V=1/60\), \(r_Y=60\), and Reid's formula gives
\(P_8-P_7=P_{16}-P_{15}=2\).
Since \(r_YV=1\), \hyperlink{result-3.12}{Lemma 3.13}
makes \(\mathcal A=R/(s)\) a domain, and its multiplication estimate gives
\[
2=\dim\mathcal A_{16}\ge2\dim\mathcal A_8-1=3.
\]
This contradiction leaves the two displayed triples. \(\square\)

For the subsequent restrictions, let \(N_m\) be the moving divisor of
the restricted \(m\)-canonical system on the original special fibre.
The coefficientwise comparison of \hyperlink{result-2.2}{Proposition 2.2} gives
\[
mD-N_m\ge_{\mathbb Q}0,\qquad
N_m^2\le mD\cdot N_m\le m^2D^2.
\tag{5.10}\label{eq:5.34}
\]
Here \(N_m\) is nef even when isolated base points remain; the upper
bounds use effectiveness and nefness.

\hypertarget{canonical-subrings-and-seventh-sections}{%
\hypertarget{section-5.4}{}\subsection{5.3. Canonical subrings and seventh sections}\label{canonical-subrings-and-seventh-sections}}

We exclude \((P_5,P_6,P_7)=(3,4,5)\). The sixth and fourteenth map
fields coincide and have index two. Let \(\iota\) be the involution,
and put \(j(a)=(-1)^m\iota^*a\) for \(a\in R_m\). Every section
through degree fourteen is fixed by \(j\).
Choose \(x\in R_4,y\in R_5,z\in R_6\) so that
\[
E_6=\langle1,t,u,h\rangle,\qquad
t=x/s^4,\quad u=y/s^5,\quad h=z/s^6.
\]
The functions \(t,u,h\) are algebraically independent, and
\(B=\mathbb C[s,x,y,z]\) is polynomial of weights \((1,4,5,6)\),
with the same fraction field as \(R^j\).
The coordinate-factor criterion of \hyperlink{result-3.13}{Lemma 3.14}
excludes every factor \(x-as^4\) from a reduced denominator of an
invariant homogeneous section: the weight sum is sixteen, so such a
factor would produce an odd section of degree thirteen.
Corollary~\hyperlink{result-canonical-polar}{3.6} also excludes \(sw=x^2\).

\hypertarget{result-5.10}{}\ResultHeading{Proposition 5.8.}
There is no seventh-canonical section satisfying \(s^2w=xy\), after
homogeneous changes of the chosen generators.

\ProofHeading{Proof.}
Such a section gives \(tu\in E_7\). For a general outer fibre
\(t=t_0\), the restriction kernel contains \(t-t_0\) and
\((t-t_0)u\). Their moving subsystem gives \(M_7\ge F+S_1\), with
\(|S_1|\) different from the outer pencil and with parameter \(u\)
on \(F\). If this curve pencil differs from the canonical pencil,
the intersection number \(\delta\) of their connected general
members is at least two. Apply \cite[Proposition 3.6(1.2)]{ref14}
with \(m_0=4,m_1=7,m_2=5,\mu=1/4,\beta=1/5\).
The sixth system separates the outer and canonical pencils, and the numerical bound is
\[
\beta^{-1}+\frac{2m_1}{\delta}
 +\mu^{-1}\left(1-\frac2\delta\right)
 =9+\frac6\delta\le12.
\]
Thus the fourteenth map would be birational. Hence \(u\) is constant
along the general canonical curve.
Put \(r=v(t)\in E_6\). By \hyperlink{result-3.2}{Lemma 3.2}(c), on the general outer fibre it is the ratio
of the independent canonical forms \(\alpha|_F\) and \((s/dt)|_F\).
Consequently \(\operatorname{trdeg}\mathbb C(t,u,r)=2\).
Writing \(r=a+bt+cu+dh\), algebraic independence of \(t,u,h\)
forces \(d=0\). This gives \(r\in E_5\), contrary to
Corollary~\hyperlink{result-canonical-polar}{3.6}. \(\square\)

Let \(\mathcal G_m\) be the saturation, in
\(f_*\omega_W^{\otimes m}\), of the subsheaf generated by all global
sections, and write \(\rho_m=\operatorname{rank}\mathcal G_m\).
Every splitting summand has nonnegative degree, since
\(\mathcal G_m\) is generically generated. Hence
\(\deg\mathcal G_m=P_m-\rho_m\), and the known spaces give
\[
\mathcal G_4=\mathcal O(1),\qquad
\mathcal G_5=\mathcal O(1)\oplus\mathcal O,\qquad
\mathcal G_6=\mathcal O(1)\oplus\mathcal O^{\oplus2}.
\tag{5.11}\label{eq:5.42}
\]

\hypertarget{result-5.11}{}\ResultHeading{Proposition 5.9.}
One has \(\rho_7=4\) and
\(\mathcal G_7=\mathcal O(1)\oplus\mathcal O^{\oplus3}\).

\ProofHeading{Proof.}
Since the outer pair is dependent over \(\mathbb C(t)\),
\(3\le\rho_7\le4\). Suppose \(\rho_7=3\). Multiplication by \(s\)
then injects \(\mathcal G_6\) into \(\mathcal G_7\) with cokernel
of length one. After normalization by \(s^m\), this is the identity
at the generic base point. An extra normalized function has form
\(P(t)+Q(t)u+R(t)h\). A modification supported at a finite
\(t=a\) gives a genuine pole in one coefficient, and therefore a
factor \(x-as^4\) in its reduced homogeneous denominator, which
has already been excluded.

The modification is supported at infinity. Frames for
\(\mathcal G_6\) are \(1,u,h\) on the finite chart and \(t,u,h\)
at infinity, so the extra function, modulo \(E_6\), has form
\[
at^2+btu+cth.
\]
If \(c\ne0\), an allowed degree-six change replaces this by \(th\),
giving the hypersurface \(xz-s^3w=0\). It has form \(UV=g\),
hence is canonical by the argument in the proof of
\hyperlink{result-3.13}{Lemma 3.14}(c). Its canonical generator has
degree \(23-10=13\), so part (b) gives an odd section of degree fourteen.
If \(c=0\) and \(b\ne0\), changing \(u\) by a multiple of \(t\)
gives \(tu\), excluded by \hyperlink{result-5.10}{Proposition 5.8}.
The remaining possibility \(t^2\) gives \(sw=x^2\), excluded by
Corollary~\hyperlink{result-canonical-polar}{3.6}.
Thus \(\rho_7=4\). The nonnegative splitting degrees, their sum one,
and \(P_7=5\) give the asserted splitting. \(\square\)

\hypertarget{result-5.12}{}\ResultHeading{Lemma 5.10.}
If \(\xi\ge2/5\), then \(D^2\le1/10\).

\ProofHeading{Proof.}
The restricted fourteenth system contains
\(|K_F+\lceil9D\rceil|\), which is nonbirational.
If \(D^2>1/10\), Langer's criterion
\cite[Remark 3.1.2 and Theorem 4.1]{ref23} gives a covering curve
\(B\) with
\[
9D\cdot B\le\frac4{1+\sqrt{1-8/(81D^2)}}<18/5.
\]
By \hyperlink{result-5.4}{Lemma 5.4}, a covering curve of canonical
degree one is a canonical-pencil member, whose \(9D\)-degree is
at least \(18/5\). Every other covering curve has the same lower
bound from \(D\ge_{\mathbb Q}A/5\). This contradiction proves the
lemma. \(\square\)

\hypertarget{result-exclude-345}{}\ResultHeading{Proposition 5.11.}
The triple \((P_5,P_6,P_7)=(3,4,5)\) is impossible under
fourteen-canonical failure.

\ProofHeading{Proof.}
The function \(r=v(t)\) restricts to a degree-one parameter of the
canonical pencil on the general outer fibre. If \(u\) is constant
on its general canonical curve, the same independence argument as
in \hyperlink{result-5.10}{Proposition 5.8} gives \(r\in E_5\), contradicting
Corollary~\hyperlink{result-canonical-polar}{3.6}. Thus \(u|_C\) is nonconstant.
The marked pole bound and the Weierstrass semigroup give
\(5\xi\ge2\). By \hyperlink{result-5.12}{Lemma 5.10}, \(D^2\le1/10\).

Put \(N=N_7\) on the original special fibre. Its map has degree two:
its field is the restriction of the common field \(L_6=L_{14}\).
By \hyperlink{result-5.11}{Proposition 5.9} its image is a nondegenerate surface in
\(\mathbb P^3\), of degree at least two. Resolving isolated base
points therefore gives moving square at least four. Equation
\eqref{eq:5.34} yields
\[
4\le N^2\le49D^2\le49/10<5.
\]
Hence \(N^2=4\), there are no base points, and the image is a
quadric. On this same original model, \hyperlink{result-2.2}{Proposition 2.2} gives
\[
K_F\cdot N\le5D\cdot N\le35D^2\le7/2.
\]
Hodge index gives \(A\cdot N\ge2\), and \(K_F-A\) is effective.
Adjunction makes \(K_F\cdot N\) even, so
\(K_F\cdot N=A\cdot N=2\). Equality in Hodge index gives
\(N\equiv2A\). Consequently
\[
7D^2\ge D\cdot N=2A\cdot D\ge2\xi\ge4/5,
\]
contrary to \(D^2\le1/10\). \(\square\)

\hypertarget{birationality-in-degree-fourteen}{%
\needspace{10\baselineskip}\hypertarget{section-6}{}\section{6. Birationality in degree fourteen}\label{birationality-in-degree-fourteen}}

\needspace{4\baselineskip}

\hypertarget{the-remaining-quadratic-configuration}{%
\hypertarget{section-6.2}{}\subsection{6.1. The remaining quadratic configuration}\label{the-remaining-quadratic-configuration}}

We now turn to the remaining case \((P_5,P_6,P_7)=(2,4,4)\). \hyperlink{result-3.19}{Proposition 3.17} and \hyperlink{result-3.20}{Corollary 3.18} supply \(k\ne0\) and sections \(x\in R_4,r_6,z_6\in R_6\) such that
\[t=x/s^4,\qquad r=r_6/s^6=v(t),\qquad z=z_6/s^6=v^2(t)-kt^2.\]
The sixth map has degree two, and its field equals the fourteenth map field:
\[[\mathbb C(X):L_6]=2,\qquad L_6=L_{14}.\tag{6.1}\label{eq:6.21}\] Its involution negates \(s\) and \(\alpha\). At every reduced transverse component of \(\operatorname{div}_Z(s)\) the section \(x\) is a unit; on all other components of multiplicity \(\ell\) it has order at least \(\ell\). There is at least one reduced transverse component. These facts use the Euler identity on the full nonreduced canonical divisor.

\needspace{3\baselineskip}

\hypertarget{result-6.6}{}\ResultHeading{Proposition 6.1.} Under the low-plurigenus and quadratic-field assumptions above, if \(z\in\langle1,t,r\rangle_{\mathbb C}\), then \(L_{13}\ne L_6\). Thus the thirteenth canonical map, and in particular the fourteenth canonical map, is birational.

\ProofHeading{Proof.} The equality \(z\in\langle1,t,r\rangle\) gives constants \(A,B,C\) with

\[
v^2t=kt^2+A+Bt+Cvt.
\tag{6.2}\label{eq:6.22}
\]

Choose \(u=u_6/s^6\) complementary to \(1,t,r\) in \(R_6/s^6\). The functions \(t,r,u\) are algebraically independent, since \(L_6\) has transcendence degree three. Define

\[
J=\frac{dt\wedge dr\wedge du}{s}\ne0.
\tag{6.3}\label{eq:6.23}
\]

It is anti-invariant. We prove \(s^{13}J\in R_{13}\) by checking prime divisors on \(Z\).

First suppose \(\alpha|_T=0\). Locally \(s=x^a h\nu\), and \(v=x^{-a}w\), where \(w\) is regular and tangent to \(T\) (the unit \(h\) is absorbed into \(w\)). Thus applying \(v\) increases a pole order by at most \(a\). If \(t\) has pole order \(n>2a\), the term \(kt^2\) in \eqref{eq:6.22} has strictly greater pole order than \(v^2t\), \(Cvt\), \(Bt\) and \(A\), which is impossible. Hence \(t\) and \(r\) have pole orders at most \(2a\) and \(3a\), respectively. The pole order of \(u\) is at most \(6a\). Alternation of three differentials loses at most one additional normal order, and division by \(s\) loses \(a\) more. Therefore

\[
\operatorname{ord}_T J\ge -(2a+3a+6a+a+1)
=-(12a+1)\ge-13a.
\tag{6.4}\label{eq:6.24}
\]

At a transverse component use the actual time \(\tau=\rho^{a+1}\) of \hyperlink{result-3.6}{Lemma 3.7}. The existing meromorphic expansion has exponents in \((a+1)^{-1}\mathbb Z\).

If \(t\) has a pole, comparison of the lowest exponents in \eqref{eq:6.22} forces its leading term to be \((6/k)\tau^{-2}\). The coefficient recursion at a shift \(q>0\) has factor

\[
(q-2)(q-3)-12=(q+1)(q-6).
\tag{6.5}\label{eq:6.25}
\]

It is nonzero for \(0<q<6\). All coefficients before the shift \(6\) are therefore constants in \(y,z_0\). In particular

\[
\operatorname{ord}_{\tau}t_y,\operatorname{ord}_{\tau}t_{z_0}\ge4,
\qquad
\operatorname{ord}_{\tau}r_y,\operatorname{ord}_{\tau}r_{z_0}\ge3.
\]

Put \(h=a+1\). Since \(u\) has \(\rho\)-pole order at most \(6a\), its transverse derivatives have \(\tau\)-order at least \(-6a/h\), and \(u_\tau\) has order at least \(-6a/h-1\). Expanding the determinant in \eqref{eq:6.23} in its time column gives

\[
\operatorname{ord}_{\tau}J\ge -6a/h,
\qquad \operatorname{ord}_T J\ge-6a.
\tag{6.6}\label{eq:6.26}
\]

If \(t\) has no pole, its Puiseux series in \(\tau\) has only integral nonnegative exponents. Indeed, the least nonintegral exponent \(p>0\) would give the uncancelled term \(p(p-1)a_p\tau^{p-2}\) on the left of \eqref{eq:6.22}. Thus \(t,r\) and their time/transverse derivatives are holomorphic in \(\tau\), and

\[
\operatorname{ord}_T J\ge-6a-(a+1)=-7a-1\ge-13a.
\tag{6.7}\label{eq:6.27}
\]

Outside \(\operatorname{div}(s)\) there is no pole. The codimension-one criterion for the reflexive sheaf \(\omega_Z^{[13]}\) now proves \(s^{13}J\in R_{13}\). This contradicts \(L_{13}=L_6\). Since \(L_{13}\) is intermediate in a quadratic extension, \(L_{13}=\mathbb C(W)\). \(\square\)

\textbf{Remark.} The stronger degree-thirteen conclusion is conditional on the sixth-canonical map already having degree two and on its natural sign in \eqref{eq:6.21}. It is not a global thirteenth-canonical theorem.

We may therefore assume that \(t,r,z\) are algebraically independent and that
\[E_6=E_7=\langle1,t,r,z\rangle,\quad t'=r,\quad r'=z+kt^2,\quad s^8z'\in R_8.\tag{6.8}\label{eq:6.28}\] The common plateau lemma removes one numerical endpoint.

\needspace{3\baselineskip}

\hypertarget{result-6.7}{}\ResultHeading{Lemma 6.2.} In this configuration \(\xi<1/2\).

\ProofHeading{Proof.} The strict curve threshold gives \(\xi\le1/2\), with equality forcing \(V=1/40\). The initial basket is
\[8(1,2)+3(1,3)+(1,4)+(2,5).\]
Its volume is \(1/20\). Enumerate its prime packings using the determinant-one rule, discarding a branch as soon as its volume is below \(1/60\). A packing strictly reduces both the number of entries and the volume, so this finite procedure is exhaustive. Among the descendants with \((P_1,\ldots,P_7)=(1,1,1,2,2,4,4)\) there is exactly one of volume \(1/40\):
\[8(1,2)+(3,8)+(3,10).\]
Here the two nontrivial packings are \((1,3)+(1,4)\rightsquigarrow(2,7)\) followed by \((1,3)+(2,7)\rightsquigarrow(3,10)\), and \((1,3)+(2,5)\rightsquigarrow(3,8)\). Completeness can also be checked from the finite budget identities \(\sum b=14\), \(\sum\gamma=5\), \(\sum e=17/120\) by listing the primitive pairs; equivalently the exact packing procedure in Appendix A gives 17 possible descendants and the unique displayed endpoint. Reid's formula gives \(P_8=P_9=6\). \hyperlink{result-3.21}{Lemma 3.19} excludes it. \(\square\)

\needspace{4\baselineskip}

\hypertarget{eighth-sections-and-the-final-jacobian}{%
\hypertarget{section-6.3}{}\subsection{6.2. Eighth sections and the final Jacobian}\label{eighth-sections-and-the-final-jacobian}}

Put \(K=\mathbb C(t)\). The functions \(r=v(t)\) and \(z\) above generate the quadratic invariant field over \(K\). The restriction parameters remain \(\mu=1/4\), \(\beta=1/5\), and degree six separates both pencils. The coordinate-factor lemma applied to \(\mathbb C[s,x,r_6,z_6]\), of weight sum seventeen, excludes every factor \(x-a s^4\) in the reduced denominator of any invariant homogeneous section: it would give an odd section of degree fourteen. The actual pole estimate on the canonical genus-two curve, and \(1/3\le\xi<1/2\), give \(2\le\lfloor8\xi\rfloor\le3\).

\needspace{3\baselineskip}

\hypertarget{result-6.8}{}\ResultHeading{Lemma 6.3.} On the generic outer fibre, \(r=v(t)\) is a degree-one parameter of the connected canonical genus-two pencil.

\ProofHeading{Proof.} On a neighbourhood of a general smooth fibre \(F\), the outer coordinate \(t\) is a local parameter on the base. The holomorphic form \(s\) has a well-defined relative restriction
\[
a=(s/dt)|_F\in H^0(F,K_F).
\]
Put \(b=\alpha|_F\). The identity \(dt\wedge\alpha=v(t)s\) gives
\[
b=r|_F\,a.
\tag{6.9}\label{eq:6.29}
\]
The forms \(a,b\) are independent: otherwise \(r\) would be constant on the generic fibre and algebraic over \(\mathbb C(t)\), contrary to \eqref{eq:6.28} and connectedness of \(f\). Since \(p_g(F)=2\), they span the canonical sections. Removing the common exceptional fixed part and descending to \(S\) preserves their ratio. Hence \(r=b/a\) is a fractional-linear coordinate on the base of the canonical pencil. \(\square\)

Write \(K=\mathbb C(t)\). The generic inner curve is thus defined over \(K(r)\). By \(1/3\le\xi<1/2\),
\[
2\le\lfloor8\xi\rfloor\le3.
\]
At a Weierstrass point on a genus-two curve, both \(H^0(2p)\) and \(H^0(3p)\) have dimension two. The functions \(1,z\) are independent over \(K(r)\) by \eqref{eq:6.28}, and belong to the space in \eqref{eq:2.10}, because \(E_6\subset E_8\). Therefore they form its basis. Every \(h\in E_8\) has a unique expression
\[
h=A_h(t,r)z+B_h(t,r),\qquad
A_h,B_h\in K(r).
\tag{6.10}\label{eq:6.30}
\]
Because we work over the generic inner curve, \eqref{eq:6.30} is an identity in \(\mathbb C(X)\). We can therefore use the same coefficients when varying the complex fibres.

We shall use the following direct consequence of the restriction criterion of Chen--Hu--Penegini. We record the numerical substitution so that its application is explicit.

\needspace{3\baselineskip}

\hypertarget{result-6.9}{}\ResultHeading{Lemma 6.4.} If the restricted eighth moving system contains
\(C+C_1\), with \(C_1\) moving and \(CC_1=\delta\ge2\), then \(\varphi_{14}\)
is birational.

\ProofHeading{Proof.} We apply \cite[Proposition 3.7]{ref14} with \(m_1=8\), \(j=1\), \(m_2=5\),
\(\mu=1/4\) and \(\beta=1/5\). For \(\delta=2\) its required integer is
\(n\ge13\); for \(\delta>2\) the bound is
\(n>9+6/\delta\le11\). Thus \(n=13\) works in both cases. \(\square\)

\needspace{3\baselineskip}

\hypertarget{result-6.10}{}\ResultHeading{Proposition 6.5 (constant transverse coefficient).} If \(\varphi_{14}\) is not birational, every coefficient \(A_h\) in \eqref{eq:6.30} is independent of \(r\):
\[
A_h(t,r)=A_h(t)\in K.
\tag{6.11}\label{eq:6.31}
\]

\ProofHeading{Proof.} Suppose that \(A_h\) depends on \(r\). Restrict to a general complex outer fibre \(F\) on which that dependence remains nonconstant and all expressions are defined. For a general constant \(c\), let \(C_c\) be the canonical member \(r=c\). The two normalized functions
\[
r-c,\qquad
h-A_h(t,c)z-B_h(t,c)
\tag{6.12}\label{eq:6.32}
\]
represent sections of the restricted eighth-canonical system. Here \(t\) is fixed on \(F\), so the displayed coefficients are complex scalars. Both sections vanish along \(C_c\), since \eqref{eq:6.30} restricts to that identity there. The normalizing section \(s^8|_F\) does not vanish identically on a general \(C_c\).

Their ratio equals
\[
\frac{A_h(t,r)-A_h(t,c)}{r-c}\,z
+\frac{B_h(t,r)-B_h(t,c)}{r-c}.
\tag{6.13}\label{eq:6.33}
\]
It does not belong to \(\mathbb C(r)\) on the general \(F\), since the coefficient of \(z\) is nonzero and \(z\) is transcendental over the inner base field. Thus the pencil generated by \eqref{eq:6.32}, after removing \(C_c\) and its remaining fixed curves, is not composed with the canonical pencil.

Resolve that pencil. Its Stein base is rational because \(q(F)=0\). If its moving divisor is a multiple of a connected general fibre, retain one such fibre \(C_1\). Then \(M_8|_F\ge C+C_1\). The map from a general canonical genus-two curve to this rational base is nonconstant, so its degree is at least two. Equivalently \(C\cdot C_1\ge2\).

\hyperlink{result-6.9}{Lemma 6.4} gives birationality in degree fourteen, a contradiction. Hence \(A_h\) is independent of \(r\). The unique coefficient in \eqref{eq:6.30} lies in \(K(r)\), so independence of \(r\) gives membership in \(K\), without a field-extension ambiguity. \(\square\)

\needspace{3\baselineskip}

\hypertarget{result-6.11}{}\ResultHeading{Proposition 6.6.} Under nonbirationality, \eqref{eq:6.30} has the more precise form
\[
h=A_h(t)z+B_h(t,r),\qquad B_h\in K[r].
\tag{6.14}\label{eq:6.34}
\]

\ProofHeading{Proof.} On the minimal surface of the generic outer fibre, all normalized sections \(h\in E_8\) have poles only on the canonical divisor \((a=0)\) selected in \hyperlink{result-6.8}{Lemma 6.3}. Indeed, before descending, their denominators are powers of \(s|_F\), and holomorphic pluricanonical sections acquire the same exceptional fixed factors under a surface blowup. The same assertion holds for \(z\in E_6\).

Since \(A_h(t)\) is constant on this surface, the rational function
\[
B_h=h-A_h(t)z\in K(r)
\]
also has no divisorial poles outside \((a=0)\). After extension to an algebraic closure of \(K\), a finite pole \(r=c\) of \(B_h\) would give a polar component of the canonical member \((b-ca=0)\). This member and \((a=0)\) share no curve, because the canonical pencil of \(S\) has no fixed curve. Such a pole is impossible. Therefore \(B_h\) has no finite pole on the rational \(r\)-line, and belongs to \(K[r]\). Polynomiality descends to \(K\) by uniqueness. \(\square\)

Apply \eqref{eq:6.34} to the particular eighth-canonical function \(w=v(z)\) from \eqref{eq:6.28}. There are rational functions of \(t\) such that
\[
w=A(t)z+\sum_{j=0}^{d}B_j(t)r^j,
\qquad A,B_j\in\mathbb C(t).
\tag{6.15}\label{eq:6.35}
\]

\needspace{3\baselineskip}

\hypertarget{result-6.12}{}\ResultHeading{Proposition 6.7.} All coefficients in \eqref{eq:6.35} are polynomials in \(t\).

\ProofHeading{Proof.} Suppose some coefficient has a finite pole at \(t=a\). Take the largest pole order among the coefficients. Since \(r,z\) are algebraically independent over \(\mathbb C(t)\), the nonzero leading polar coefficients cannot cancel between the distinct monomials
\[
z,1,r,r^2,\ldots,r^d.
\]
Homogenize \eqref{eq:6.35} using \(t=x/s^4\), \(r=r_6/s^6\) and \(z=z_6/s^6\), and multiply by \(s^8\). In the reduced homogeneous fraction for \(w_8\) over \(B\), a factor \(x-a s^4\) remains in the denominator. One can check noncancellation after inverting \(s\): the quotient by \(x-a s^4\) still has independent variables \(r_6,z_6\).

This is forbidden by the coordinate-factor exclusion, which applies to invariant sections in arbitrary degree and in particular to \(w_8\in R_8\). Thus no coefficient has a finite pole. A rational function of one variable over \(\mathbb C\) with no finite pole is polynomial. \(\square\)

We have therefore proved
\[
v(z)=A(t)z+P(t,r),\qquad A\in\mathbb C[t],\quad P\in\mathbb C[t,r].
\tag{6.16}\label{eq:6.36}
\]

Define
\[
J=\frac{dt\wedge dr\wedge dz}{s}
=\frac{dt\wedge d(vt)\wedge d(v^2t)}{s}.
\tag{6.17}\label{eq:6.37}
\]
The second equality follows because \(dz=d(v^2t)-2kt\,dt\). By \eqref{eq:6.28}, the function \(J\) is nonzero, and \(\iota(J)=-J\). It remains to prove that \(s^{14}J\) is a global section. We do this by checking its order along each prime divisor, beginning with the components at which the estimate is automatic.

\hypertarget{result-6.13}{}\ResultHeading{Proposition 6.8 (the critical-pole criterion).}
At a simple transverse component, define the obstruction
\(da_2\wedge da_3\) using the expansion \eqref{eq:6.41}.
At a double transverse component where \(t\) has pole order six,
define the obstruction \(db_0\wedge db_1\) using \eqref{eq:6.45}.
The rational section \(s^{14}J\) extends to \(R_{14}\) if and only if
all these obstructions vanish. Every other prime divisor satisfies
the required order automatically.

\ProofHeading{Proof.}
For functions with signed orders \(a_1,a_2,a_3\) along \(T\), alternation gives
\[
\operatorname{ord}_T\frac{df_1\wedge df_2\wedge df_3}{s}
\ge a_1+a_2+a_3-\ell-1.
\tag{6.18}\label{eq:6.38}
\]
At most one differential contributes an extra normal order.

At a tangent component, the componentwise order bounds above bound the poles of \(t,vt,v^2t\) by \(3\ell,4\ell,5\ell\). Therefore
\[
\operatorname{ord}_T(J)\ge-13\ell-1\ge-14\ell.
\tag{6.19}\label{eq:6.39}
\]

At a transverse component of multiplicity \(\ell\ge2\), the pole order \(n\) of \(t\) in \(\rho\) satisfies \(n\le3\ell\); use \(n=0\) if \(t\) is regular. In the actual time \(\tau=\rho^{\ell+1}\), expand the already meromorphic function in a convergent Puiseux series. In a nonzero determinant term, the two transverse columns cannot use the same exponent: their coefficient vectors would be proportional. The next exponent differs by at least \(1/(\ell+1)\). The total derivative loss in \eqref{eq:6.37} is four time derivatives. Thus
\[
\operatorname{ord}_T(J)
\ge-3n-4(\ell+1)+1.
\tag{6.20}\label{eq:6.40}
\]
For \(\ell\ge3\), this is at least \(-13\ell-3\ge-14\ell\). For \(\ell=2\) and \(n\le5\), it is at least \(-26>-28\). Only a simple transverse component, or a double transverse component with \(n=6\), can obstruct degree fourteen.

Let \(\ell=1\) and \(\tau=\rho^2\). The unit condition for \(x\) gives a leading pole \(t=a_0\tau^{-2}+\cdots\). The condition that \(z=t''-kt^2\) have \(\rho\)-pole at most six forces
\[
a_0=\frac6k,\qquad
t=a_0\tau^{-2}+a_2\tau^{-1}
+a_3\tau^{-1/2}+O(1).
\tag{6.21}\label{eq:6.41}
\]
Indeed, cancellation of the \(\tau^{-4}\) term fixes \(a_0\), and the coefficient at \(\tau^{-7/2}\) then forces the missing \(\tau^{-3/2}\) coefficient of \(t\) to vanish. The functions \(a_2,a_3\) are holomorphic in the transverse variables at a general point.

For a time-column exponent \(p\) and two transverse exponents \(q,r_0\), the numerical determinant is
\[
\det\begin{pmatrix}
p&1&1\\
p(p-1)&q&r_0\\
p(p-1)(p-2)&q(q-1)&r_0(r_0-1)
\end{pmatrix}
=p(q-p+1)(r_0-p+1)(r_0-q).
\tag{6.22}\label{eq:6.42}
\]
At \((p,q,r_0)=(-2,-1,-1/2)\) its value is \(-5\). Therefore
\[
[\tau^{-15/2}]J
=-\frac{30}{k}\,
\frac{da_2\wedge da_3}{du_1\wedge du_2}.
\tag{6.23}\label{eq:6.43}
\]
There is no earlier term, and subsequent exponents increase in steps of \(1/2\). Consequently
\[
\operatorname{ord}_T(J)\ge-14
\quad\Longleftrightarrow\quad da_2\wedge da_3=0.
\tag{6.24}\label{eq:6.44}
\]
We call the component simple critical when this wedge is nonzero.

Let \(\ell=2\), \(n=6\) and \(\tau=\rho^3\). Write
\[
t=b_0\tau^{-2}+b_1\tau^{-5/3}
+O(\tau^{-4/3}),\qquad b_0\ne0.
\tag{6.25}\label{eq:6.45}
\]
Using \eqref{eq:6.42} with \((-2,-2,-5/3)\) gives
\[
[\tau^{-29/3}]J
=-\frac89 b_0\,
\frac{db_0\wedge db_1}{du_1\wedge du_2}.
\tag{6.26}\label{eq:6.46}
\]
The next possible order in \(\rho\) is at least \(-28\). Hence
\[
\operatorname{ord}_T(J)\ge-28
\quad\Longleftrightarrow\quad db_0\wedge db_1=0.
\tag{6.27}\label{eq:6.47}
\]
We call the component double critical when this wedge is nonzero.

The two displayed wedge conditions are therefore precisely the remaining extension obstructions.

There are no poles away from \(\operatorname{div}_Z(s)\): all normalized functions in \eqref{eq:6.37} are regular there in codimension one. Since \(Z\) is normal and the pluricanonical sheaf is reflexive, the codimension-one test suffices. No assertion about a smooth index-one cover of a basket point is needed. \(\square\)

\needspace{3\baselineskip}

\hypertarget{result-6.14}{}\ResultHeading{Proposition 6.9.} Equation \eqref{eq:6.36} excludes both types of critical component.

\ProofHeading{Proof.} The vector field in coordinates \((t,r,z)\) is
\[
v=r\partial_t+(z+kt^2)\partial_r+(A(t)z+P(t,r))\partial_z.
\]
Its divergence is \(A(t)\). At either critical component \(t\) has time-pole
two, so \hyperlink{result-3.9}{Lemma 3.10}(a) gives \(A\in\mathbb C\).

At a simple critical component, use \eqref{eq:6.41}, with \(a_0=6/k\), and put
\(H=r^2-(2k/3)t^3\). The leading terms are
\[
r=-2a_0\tau^{-3}+\cdots,\qquad
H=-8a_0a_2\tau^{-5}+\cdots,\qquad
z'=30a_2\tau^{-4}+\cdots.
\]
Every \(P\in\mathbb C[t,r]\) is uniquely a sum of \(t^iH^jr^\varepsilon\),
with \(\varepsilon=0,1\). At fixed weight \(2i+5j+3\varepsilon\), the leading
coefficients are nonzero constants times distinct powers of \(a_2\).
Since \(a_2\) is nonconstant, they cannot cancel. Now \(Az\) has pole at most
three and \(z'\) has pole four, so
\(P\in\langle1,t,r,t^2\rangle\). The coefficient \(e\) of \(t^2\) would satisfy
\(30a_2=e a_0^2\), contradicting nonconstancy.

At a double critical component, \eqref{eq:6.45} gives
\[
r=-2b_0\tau^{-3}+\cdots,\quad
z=b_0(6-kb_0)\tau^{-4}+\cdots,\quad
z'=-4b_0(6-kb_0)\tau^{-5}+\cdots.
\]
At fixed weight \(2i+3j\), the leading coefficients of the monomials
\(t^ir^j\) are distinct powers \((-2)^jb_0^{i+j}\); they are independent
because \(b_0\) is nonconstant. Hence
\(P\in\langle1,t,r,t^2,tr\rangle\). If \(e\) is the coefficient of \(tr\),
comparison at pole five yields
\(-24b_0+(4k+2e)b_0^2=0\), again impossible. \(\square\)

\ProofHeading{Proof of \hyperlink{result-1.1}{Theorem 1.1}.}
Suppose that \(\varphi_{14}\) is not birational.
By \hyperlink{result-4.4}{Proposition 4.4} and the exclusions for
fibres with canonical square two, the four-canonical pencil has
\((1,2)\) surfaces over \(\mathbb P^1\) as general fibres.
The analysis of the sixth system and
\hyperlink{result-exclude-345}{Proposition 5.11} leave only
\((P_5,P_6,P_7)=(2,4,4)\).

\hyperlink{result-6.6}{Proposition 6.1} reduces us to the independent
coordinates \(t,r,z\) in \eqref{eq:6.28}.
Propositions~6.5--6.7 give the polynomial equation
\eqref{eq:6.36}. By \hyperlink{result-6.14}{Proposition 6.9}, neither
type of critical component occurs. Thus
\hyperlink{result-6.13}{Proposition 6.8} gives
\(s^{14}J\in R_{14}\).
The nonzero anti-invariant function \(J\) contradicts
\(L_{14}=L_6\).
Hence \(\varphi_{14}\) is birational. Multiplication by
\(s^{m-14}\) embeds its normalized sections into \(E_m\), so
\(\varphi_m\) is birational for every \(m\ge14\). \(\square\)

\needspace{4\baselineskip}

\hypertarget{sharpness}{%
\needspace{10\baselineskip}\hypertarget{section-7}{}\section{7. Sharpness}\label{sharpness}}

Let \(X_{28}\) be a general quasismooth well-formed hypersurface in \(\mathbb P(1,3,4,5,14)\), with coordinates \(s,x,y,z,w\) of the indicated weights. This is a canonical threefold of general type; see \cite{ref20} and \cite[Example 1.2(3)]{ref14}. Weighted adjunction and intersection give
\[
K_{X_{28}}\cong\mathcal O_{X_{28}}(1),\qquad
p_g(X_{28})=1,\qquad
K_{X_{28}}^3=\frac{28}{1\cdot3\cdot4\cdot5\cdot14}=\frac1{30}.
\]
After completing the square, its equation has the form \(w^2=F_{28}(s,x,y,z)\). Every section of degree at most thirteen is independent of \(w\). In particular, \(\varphi_{13}\) factors through the double cover of \(\mathbb P(1,3,4,5)\) and is not birational. In degree fourteen, the ratios of \(s^{14}\), \(s^{11}x\), \(s^{10}y\), \(s^9z\), and \(w\) generate the function field on \(s\ne0\). Thus \(\varphi_{14}\) is birational, while \(\varphi_{13}\) is not. This proves that the bound in \hyperlink{result-1.1}{Theorem 1.1} is sharp.

\needspace{4\baselineskip}

\hypertarget{appendix-a.-the-volume-one-fortieth-endpoint}{%
\section{Appendix A. The volume-one-fortieth endpoint}\label{appendix-a.-the-volume-one-fortieth-endpoint}}

In the remaining configuration, the canonical unpacking identities of \cite[Section 3]{ref9} give
\[B^{(5)}=8(1,2)+3(1,3)+(1,4)+(2,5),\qquad V(B^{(5)})=1/20.\]
Every basket is obtained from this one by prime packings. At each step replace two pairs \((b,r),(c,u)\) with \(|bu-cr|=1\) by \((b+c,r+u)\). The volume decreases by \(1/(ru(r+u))\). Stop a branch as soon as its volume is below \(1/60\), and keep its distinct nodes only once. This procedure terminates, since each step decreases the number of entries. It is exhaustive because every prime-packing chain is one of its branches and volume never increases. Finally impose \((P_1,\ldots,P_7)=(1,1,1,2,2,4,4)\) using Reid's formula. The resulting table is

\begin{longtable}[]{@{}
  >{\raggedright\arraybackslash}p{(\columnwidth - 6\tabcolsep) * \real{0.5800}}
  >{\raggedleft\arraybackslash}p{(\columnwidth - 6\tabcolsep) * \real{0.1400}}
  >{\raggedleft\arraybackslash}p{(\columnwidth - 6\tabcolsep) * \real{0.1400}}
  >{\raggedleft\arraybackslash}p{(\columnwidth - 6\tabcolsep) * \real{0.1400}}@{}}
\toprule\noalign{}
\begin{minipage}[b]{\linewidth}\raggedright
Basket
\end{minipage} & \begin{minipage}[b]{\linewidth}\raggedleft
\(V\)
\end{minipage} & \begin{minipage}[b]{\linewidth}\raggedleft
\(P_8\)
\end{minipage} & \begin{minipage}[b]{\linewidth}\raggedleft
\(P_9\)
\end{minipage} \\
\midrule\noalign{}
\endhead
\bottomrule\noalign{}
\endlastfoot
\(3(1,2)+3(1,3)+(1,4)+(7,15)\) & \(1/60\) & 6 & 7 \\
\(7(1,2)+(1,3)+(3,7)+(3,10)\) & \(2/105\) & 5 & 6 \\
\(4(1,2)+3(1,3)+(1,4)+(6,13)\) & \(1/52\) & 6 & 7 \\
\(5(1,2)+3(1,3)+(1,4)+(5,11)\) & \(1/44\) & 6 & 7 \\
\(7(1,2)+2(1,3)+(2,7)+(3,7)\) & \(1/42\) & 5 & 6 \\
\(8(1,2)+(3,8)+(3,10)\) & \(1/40\) & 6 & 6 \\
\(8(1,2)+(2,7)+(4,11)\) & \(2/77\) & 6 & 6 \\
\(6(1,2)+3(1,3)+(1,4)+(4,9)\) & \(1/36\) & 6 & 7 \\
\(8(1,2)+(1,3)+(2,7)+(3,8)\) & \(5/168\) & 6 & 6 \\
\(8(1,2)+(2,5)+(4,13)\) & \(2/65\) & 6 & 7 \\
\(8(1,2)+(1,3)+(2,5)+(3,10)\) & \(1/30\) & 6 & 7 \\
\(8(1,2)+(1,4)+(5,14)\) & \(1/28\) & 7 & 7 \\
\(7(1,2)+3(1,3)+(1,4)+(3,7)\) & \(1/28\) & 6 & 7 \\
\(8(1,2)+(1,3)+(1,4)+(4,11)\) & \(5/132\) & 7 & 7 \\
\(8(1,2)+2(1,3)+(2,5)+(2,7)\) & \(4/105\) & 6 & 7 \\
\(8(1,2)+2(1,3)+(1,4)+(3,8)\) & \(1/24\) & 7 & 7 \\
\(8(1,2)+3(1,3)+(1,4)+(2,5)\) & \(1/20\) & 7 & 8 \\
\end{longtable}

Only \(8(1,2)+(3,8)+(3,10)\) has volume \(1/40\), and there \(P_8=P_9=6\). The rows are numerical possibilities; their enumeration makes no assertion about geometric realization. The following computation uses exact rational arithmetic and verifies the list directly from the packing rule and the displayed initial basket. An executable copy is supplied as the ancillary file \texttt{verify\_endpoint.py}, together with its expected output and running instructions.

\begin{verbatim}
from fractions import Fraction as Q

def volume(B):
    return sum((Q(b*b, r) for b, r in B), Q(0)) - 6

def P(B, m):
    if m == 1:
        return Q(1)
    correction = sum((Q((j*b % r)*(r-j*b % r), 2*r)
                      for b, r in B for j in range(1, m)), Q(0))
    return Q(m*(m-1)*(2*m-1), 12)*volume(B) - (2*m-1) + correction

initial = tuple(sorted([(1,2)]*8 + [(1,3)]*3 + [(1,4),(2,5)]))
todo, seen, rows = [initial], {initial}, []
while todo:
    B = todo.pop()
    if [P(B,m) for m in range(1,8)] == [1,1,1,2,2,4,4]:
        rows.append(B)
    for i in range(len(B)):
        for j in range(i+1, len(B)):
            b,r = B[i]
            c,s = B[j]
            if abs(b*s-c*r) != 1:
                continue
            C = tuple(sorted([x for h,x in enumerate(B)
                              if h not in (i,j)] + [(b+c,r+s)]))
            if C not in seen and volume(C) >= Q(1,60):
                seen.add(C)
                todo.append(C)
assert len(rows) == 17
endpoint = [B for B in rows if volume(B) == Q(1,40)]
assert len(endpoint) == 1
assert P(endpoint[0],8) == P(endpoint[0],9) == 6
\end{verbatim}

\end{document}